%% file: tex-me.tex
\documentclass[10pt,draft]%
{amsart}
\usepackage{latexsym,amsmath,amssymb}
\usepackage{graphicx,epsf}
\input{amsart-title}
\usepackage{youngtab} 

\let\storebaselineskip=\baselineskip

\usepackage{ifthen}     

\newcommand\wantcolor{no}

\newcommand\finalized{yes}

\newcommand\ignore[1]{}
\usepackage{color}
\providecommand\wantcolor{yes}   %
\ifthenelse{\equal{\wantcolor}{no}}{\renewcommand\color[1]{}}{}
\definecolor{backgroundyellow}{cmyk}{.2,.1,.8,.2}
\definecolor{backgroundblue}{rgb}{0,0,1}
\definecolor{backgroundred}{rgb}{1,0,0}
\definecolor{backgroundmagenta}{cmyk}{0,1,0,0}
\definecolor{GoodForInverseVideo}{rgb}{.6,.8,1}
\definecolor{myyellow}{rgb}{1,1,0}

\newcommand\mysection{\section}

\newcommand\mysubsection{\subsection}

\newcommand\mysubsubsection[1]{%
\subsubsection{\sffamily\upshape\mdseries #1}}
\newcommand\mysss{\mysubsubsection}

\newtheorem{annotation}{Annotation}
\newtheorem{theorem}[annotation]{
		Theorem}
\newtheorem{lemma}[annotation]{
		Lemma}
\newtheorem{corollary}[annotation]{
		Corollary}
\newtheorem{proposition}[annotation]{
		Proposition}
\newtheorem{definition}[annotation]{
		Definition}
\newcommand\bdefinition{\begin{definition}\begin{rm}}
\newcommand\edefinition{\end{rm}\end{definition}}
\newtheorem{example}[annotation]{
		Example}
\newcommand\bexample{\begin{example}\begin{rm}}
\newcommand\eexample{\end{rm}\hfill$\Box$\end{example}}
\newtheorem{notation}[annotation]{
		Notation}
\newcommand\bnotation{\begin{notation}\begin{rm}}
\newcommand\enotation{\end{rm}\end{notation}}

\newtheorem{remark}[annotation]{
		Remark}
\newcommand\bremark{\begin{remark}
\begin{upshape}}
\newcommand\eremark{\end{upshape}
\end{remark}}
\newenvironment{myproof}{%
\vspace{-1.25\parskip}\par\noindent{\scshape Proof:}\begin{rm}}{\hfill$\Box$\end{rm}\vspace{\parskip}}

\numberwithin{equation}{subsection}
\numberwithin{figure}{subsection}

\providecommand\finalized{no}
\ifthenelse{\equal{\finalized}{yes}}{}{\marginparwidth=2cm}
\ifthenelse{\equal{\finalized}{yes}}%
{\newcommand\mylabel[1]{\label{#1}}}%
{\newcommand\mylabel[1]{\label{#1}\marginpar{[{\ttfamily\upshape\tiny #1}]}}}
\usepackage{ragged2e}
\ifthenelse{\equal{\finalized}{yes}}%
{\newcommand\checked[1]{}}%
{\newcommand\checked[1]{\marginpar{[{\ttfamily\upshape\tiny CHECKED: #1}]}}}
\ifthenelse{\equal{\finalized}{yes}}%
{\newcommand\spellchecked[1]{}}%
{\newcommand\spellchecked[1]{\marginpar{[{\ttfamily\upshape\tiny SPELLCHECKED: #1}]}}}

\providecommand\version{public}   
\ifthenelse{\equal{\version}{public}}%
{\newcommand\mcomment[1]{}}%
{\newcommand\mcomment[1]{\marginpar{{\raggedright\sffamily\upshape\tiny #1}}}}
\ifthenelse{\equal{\version}{public}}%
{\newcommand\fcomment[1]{}}%
{\newcommand\fcomment[1]{\footnote{#1}}}

\input{definitions}

\begin{document}
\input{abstract}
\maketitle
\input{introorg}

\input{comb}

\input{hecke}
\input{invapp}
\input{tabloidproof}

\input{proof}
\input{matrixg}
\input{notation}

\input{gandgram}

\input{gdet}
\input{poschar}
\setlength\parskip{0\storebaselineskip}
\input{bibliography}
\end{document}

%% file: amsart-title.tex
\title[RSK bases and Kazhdan-Lusztig cells]%
{RSK bases and Kazhdan-Lusztig cells}
\author{K.~N.~Raghavan}
	\address{Institute of Mathematical Sciences, 
                 C.~I.~T.~Campus, 
		Chennai 600\,113, INDIA}
	\email{knr@imsc.res.in}

\author{Preena Samuel}
	\address{Institute of Mathematical Sciences, 
                C.~I.~T.~Campus, 
		Chennai 600\,113, INDIA}
	\email{preena@imsc.res.in}

\author{K.~V.~Subrahmanyam}
	\address{Chennai Mathematical Institute, 
		Plot No. H1, SIPCOT IT Park\newline
		\hbox{\hspace{\parindent}}%
		Padur Post, Siruseri 603\,103, 
                Tamilnadu, INDIA}
	\email{kv@cmi.ac.in}

%% file: definitions.tex
\newcommand\field{k}
\newcommand\base{k}
\newcommand\st{\,|\,}

\newcommand\sign{\textrm{sgn}}

\DeclareMathOperator{\End}{End}
\newcommand\integers{\mathbb{Z}}
\newcommand\rational{\mathbb{Q}}
\newcommand\hash{\dag}

\newcommand\lnd{\lambda(n,d)}
\newcommand\rep\rho
\newcommand\repl{\rep_\lambda}
\newcommand\replnd{\rep_{\lnd}}

\newcommand\Hom{\textup{Hom}}

\newcommand\glv{\GL(V)}

\newcommand\EndV{\Endo{V}}
\newcommand\tensor\otimes
\newcommand\symm{{\mathfrak S}}

\newcommand\jnd{J(n,d)}
\newcommand\groupring{\fieldc\symmn}

\newcommand\fieldc{\mathbb{C}}

\newcommand\symmn{\symm_n}

\newcommand\length{\ell}
\newcommand\lcmod{L}
\newcommand\rcmod{R}
\newcommand\rcmodl{\rcmod(\lambda)}
\newcommand\rcmodlc{\rcmodl_\fieldc}
\newcommand\rcmodlp{\rcmod(\lambdap)}
\newcommand\rcmodlk{\rcmodl_\field}
\newcommand\rcmodm{\rcmod(\mu)}
\newcommand\rcmodmk{\rcmodm_\field}

\newcommand\tabloid{{\mathcal T}}
\newcommand\tabloidl{\tabloid_\lambda}
\newcommand\itabl{\integers\tabloidl}

\newcommand\dija{{\mathfrak D}_\lambda}
\newcommand\wnotl{w_{0,\lambda}}
\newcommand\wnotlp{w_{0,\lambdap}}
\newcommand\vwnotlp{v_{\wnotlp}}
\newcommand\Wl{W_\lambda}

\newcommand\xl{x_\lambda}
\newcommand\xlp{x_{\lambdap}}
\newcommand\yl{y_\lambda}
\newcommand\ylp{y_{\lambdap}}
\newcommand\lambdap{\lambda'}
\newcommand\wl{w_\lambda}
\newcommand\vwl{v_{\wl}}

\newcommand\wlp{w_{\lambdap}}
\newcommand\zl{z_\lambda}
\newcommand\twl{T_{\wl}}
\newcommand\twlp{T_{\wlp}}

\newcommand\perml{M^\lambda}
\newcommand\permlp{M^{\lambdap}}
\newcommand\specht{S}
\newcommand\spechtm{\specht^\mu}
\newcommand\spechtl{\specht^\lambda}
\newcommand\spechtlp{\specht^{\lambdap}}
\newcommand\spechtlk{\spechtl_k}
\newcommand\tupl{t^\lambda}
\newcommand\tlowl{t_\lambda}
\newcommand\tlowlp{t_{\lambdap}}

\DeclareMathOperator{\GL}{GL}
\DeclareMathOperator{\Endo}{End}
\DeclareMathOperator\klla{{\leftarrow}_L}
\DeclareMathOperator{\kll}{{\leq_L}}
\DeclareMathOperator\kllt{{<_L}}

\DeclareMathOperator\klrlt{{<_R}}
\DeclareMathOperator\kllsim{{\sim_L}}

\DeclareMathOperator\klr{{\leq_R}}
\DeclareMathOperator\klrsim{{\sim_R}}
\DeclareMathOperator\kllr{{\leq_{LR}}}
\DeclareMathOperator\kllrsim{{\sim_{LR}}}
\DeclareMathOperator\klllneq{{\lneq_{L}}}
\DeclareMathOperator\klrlneq{{\lneq_{R}}}
\DeclareMathOperator\kllrlneq{{\lneq_{LR}}}
\DeclareMathOperator\krsshape{KRS-shape}
\newcommand\domleq{\trianglelefteq}
\newcommand\doml{\triangleleft}
\newcommand\domr{\triangleright}
\newcommand\domgeq{\trianglerighteq}

\newcommand\hecke{\mathcal H}
\newcommand\heckek{\hecke_{\field}}
\newcommand\shape{RSK-shape}
\DeclareMathOperator\rskshape{\textup{RSK-shape}}
\newcommand\partition{\vdash}

\renewcommand\bar\overline

\newcommand\subsets{\mathcal{S}}
\newcommand\subsett{\mathcal{T}}
\newcommand\matrixg{\mathcal{G}}
\newcommand\matrixgl{\matrixg(\lambda)}
\newcommand\rhol{\rho_\lambda}
\newcommand\dl{d(\lambda)}

\newcommand\basise{\mathfrak{e}}
\newcommand\bigmatrix{\mathbb{G}}
\newcommand\bigmatrixl{\bigmatrix(\lambda)}
\newcommand\bigmatrixlk{\bigmatrix(\lambda)_\base}
\newcommand\form{\langle\ ,\ \rangle}

\newcommand\ring{k}

\newcommand\lambdad{\zeta(d)}
\newcommand\wnotlambdad{w_{0,\lambdad}}
\newcommand\cwnotld{C_{\wnotlambdad}}
\newcommand\cwnotl{C_{\wnotl}}
\newcommand\wnotld{w_{0,\lambdad}}

\newcommand\lambdadp{\lambdad^\prime}

\newcommand\gramdet{\det(\lambda)}
\newcommand\qintv[1]{[#1]_v}
\newcommand\qintq[1]{[#1]_q}
\newcommand\qintqfac[1]{[#1]_q^!}

\newcommand\ideali{\mathfrak{I}}

\newcommand\nlambda{N^{\lambda}}
\newcommand\nhatl{\hat{N}^{\lambda}}
\newcommand\sblp{\tilde{\specht}^{\lambdap}}
\newcommand\sbl{\tilde{\specht}^{\lambda}}
\newcommand\forml{{\form}_\lambda}
\newcommand\dual{\textup{dual}}

\newcommand\spechtlpk{\spechtlp_\ring}

%% file: abstract.tex
\begin{abstract} 
From the combinatorial
characterizations of the right, left, and two-sided Kazhdan-Lusztig
cells of the symmetric group,   `RSK bases' are constructed for 
certain quotients by two-sided ideals of the group ring and the Hecke algebra.
Applications to invariant theory, over various base rings, 
of the general linear group and
representation theory, both ordinary and modular, 
of the symmetric group are discussed.
\end{abstract}

%% file: introorg.tex
\mysection{Introduction: summary and organization of results}%
\mylabel{s:intro}
\noindent
The starting point of the work described in this paper is a question in
classical invariant theory (\S\ref{ss:iinv}).
It leads naturally to questions 
about representations of the symmetric group over the complex numbers
(\S\ref{ss:itabloid},~\S\ref{ss:iirrep}) and over algebraically closed
fields of positive characteristic (\S\ref{ss:iposchar}),  and 
in turn to the computation of the determinant of a certain matrix
encoding the multiplication of Kazhdan-Lusztig basis elements of
the Hecke algebra (\S\ref{ss:icomputation}), 
using which one can recover a well-known criterion
for the irreducibility of Specht modules over fields of positive
characteristic~(\S\ref{ss:icarter}).

For the sake of readability,   we have tried, to the extent possible,
to keep the proofs of our results independent of each other. 
So sections~\ref{s:invapp}--\ref{s:proof} can be read without
reference to one another.    

RSK stands for Robinson-Schensted-Knuth.   
\mysubsection{Motivation from invariant theory}\mylabel{ss:iinv}
We begin by recalling a basic theorem of classical invariant 
theory. 
Let $\base$ be a commutative ring with identity and 
$V$ a free $\base$-module of finite rank $d$.
Let $\glv$ denote the group of $\base$-automorphisms of~$V$, and consider
the diagonal action of $\glv$ on $V^{\tensor n}$.   Let $\symm_n$ denote
the symmetric group of bijections of the set $\{1,\ldots,n\}$ 
and $\base\symm_n$ the group ring of~$\symmn$
with coefficients in~$\base$.      There is a natural action 
of~$\symm_n$ on~$V^{\tensor n}$ by permuting the factors:  more precisely,
$(v_1\tensor\cdots\tensor v_n)\cdot\sigma:=v_{1\sigma}\tensor\cdots
\tensor v_{n\sigma}$ (all actions are on the right by convention).
This action commutes with the action of~$\glv$,
and so the $\base$-algebra map $\Theta_n:
\base\symm_n\to\Endo_k{V^{\tensor n}}$ defining
the action of~$\symm_n$ has image in the space 
$\Endo_{\glv}{V^{\tensor n}}$ 
of $\glv$-endomorphisms of~$V^{\tensor n}$.   

We have the following result (see~\cite[Theorems~4.1,~4.2]{dp}):
\begin{quote}
Assume the following: if $f(X)$ is an element of degree~$n$ 
of the polynomial ring $\base[X]$ in one variable over~$\base$
that vanishes as a function on~$\base$,  then $f(X)$ is identically zero.
(This holds for example  when $\base$ is an infinite 
field, no matter what~$n$ is.)
Then the $\base$-algebra homomorphism $\Theta_n$ maps 
onto~$\Endo_{\glv}{V^{\tensor n}}$ and
its kernel is the two-sided ideal~$J(n,d)$ defined as follows:
\begin{itemize}
\item $J(n,d):=0$ if $d\geq n$;
\item if $d<n$,  then it is the two-sided ideal
generated
by the element $y_d:=\sum_{\tau\in\symm_{d+1}}(\sign\,\tau)\tau$,  
where $\symm_{d+1}$
is the subgroup of~$\symm_n$ consisting of the permutations that fix point-wise
the elements $d+2$, \ldots, $n$, and $\sign\,\tau$ denotes the sign of~$\tau$.%
\footnote{
The subgroup~$\symm_{d+1}$ could be taken to be that consisting
of the permutations that fix point-wise any arbitrarily fixed set of
$n-d-1$ elements.}
\end{itemize}
\end{quote}
Thus~$\base\symm_n/J(n,d)$ gets identified with the algebra of
$\glv$-endomorphisms of~$V^{\tensor n}$ (under the mild assumption on~$\base$ mentioned above),  and it is of invariant theoretic
interest to ask:
\begin{quote}
Is there a natural choice of a $\base$-basis for~$\base\symm_n/J(n,d)$?
\end{quote}

Our answer: 
\begin{theorem}\mylabel{t:invariant}
Let $\base$ be any commutative ring with identity.
  Those permutations $\sigma$ of~$\symm_n$ such that the sequence
$1\sigma$, \ldots, $n\sigma$ has no decreasing sub-sequence of
length more than~$d$ form a basis for $\base\symm_n/J(n,d)$.
\end{theorem}
\noindent
The proof of the theorem will be given in~\S\ref{s:invapp}.   It involves
the Hecke algebra of the symmetric group and its Kazhdan-Lusztig basis.   
Some further comments on the proof can be found in~\S\ref{ss:ipf}.

The theorem enables us to:
\begin{itemize}
\item obtain a 
$\base$-basis, closed under multiplication, 
for the subring of $\glv$-invariants of the tensor algebra
of~$V$~(\S\ref{ss:tensor}).
\item when $\base$ is a field of characteristic~$0$,  to limit the permutations
  in the well-known description 
(\cite{procesi},~\cite{razmyslov}) 
of a spanning set for polynomial $\glv$-invariants 
of several matrices~(\S\ref{ss:tensor});
or, more generally,
to limit the permutations
  in the description in~\cite{dks} 
of a spanning set by means of `picture invariants' for polynomial 
$\glv$-invariants of several tensors~(\S\ref{ss:poly}).
\end{itemize}
\mysubsection{A question about tabloid representations}\mylabel{ss:itabloid}
\noindent
Let us take the base ring~$\base$ in~\S\ref{ss:iinv} to be the
field~$\fieldc$ of complex numbers.    Then the ideal~$J(n,d)$ has a
representation theoretic realization as we now briefly recall
(see~\S\ref{sss:jndkerreplnd} 
for the justification).   
Let $\lnd$ be the unique partition of~$n$
with at most~$d$ parts that is smallest in the dominance order 
(\S\ref{sss:lnd}).   Consider the 
linear representation of~$\symm_n$
on the free vector space $\fieldc\tabloid_{\lnd}$ generated by
tabloids of shape $\lnd$ (\S\ref{ss:tabloid}).   The ideal~$J(n,d)$
is the kernel of the $\fieldc$-algebra map $\fieldc\symm_n\to
\Endo_{\fieldc}\fieldc\tabloid_{\lnd}$ defining this representation.

Replacing the special partition~$\lnd$ above by an arbitrary one~$\lambda$
of~$n$ (\S\ref{ss:partshape}) and considering the $\fieldc$-algebra map
$\repl:\fieldc\symm_n\to\Endo_\fieldc{\fieldc\tabloid_\lambda}$ defining the
linear representation of~$\symm_n$ on the space~$\fieldc\tabloid_\lambda$ 
generated by tabloids of shape~$\lambda$,
we ask:
\begin{quote}
  Is there a natural set of permutations that form a $\fieldc$-basis
for the group ring $\fieldc\symm_n$ modulo the kernel of the map~$\repl$?
Equivalently, one could demand that the images of the permutations
under~$\repl$ form a basis for the image.
\end{quote}

Our answer:
\begin{theorem}
  \mylabel{t:tabloid}
Permutations of~\shape~$\mu$, as $\mu$ varies over partitions that
dominate~$\lambda$, form a
$\fieldc$-basis of $\fieldc\symm_n$ modulo the kernel of 
$\repl:\fieldc\symm_n\to
\Endo{\fieldc\tabloidl}$.   
\end{theorem}
\noindent 
The dominance order on partitions is the usual one (\S\ref{ss:domination}).
The \shape\ of a permutation is defined in terms of the RSK-correspondence
(\S\ref{ss:rsk}).  As follows readily from the definitions, 
the shape of a permutation~$\sigma$ dominates the 
partition~$\lnd$ precisely when
$1\sigma$, \ldots, $n\sigma$ has no decreasing sub-sequence of 
length exceeding~$d$.    Thus, in the case when the base ring is the
complex field, Theorem~\ref{t:invariant} follows from Theorem~\ref{t:tabloid}.

The proof of the theorem will be given in~\S\ref{s:tabloidproof}.
Like that of Theorem~\ref{t:invariant}, it too involves the Hecke 
algebra of the symmetric group and its Kazhdan-Lusztig basis.
Some further comments on the proof can be found in~\S\ref{ss:ipf}.

The theorem holds also over the integers and over fields of 
characteristic~$0$---as can
be deduced easily from the complex case (see~\S\ref{ss:tabloidz})---but it is not true in general
over a field of positive characteristic: see Example~\ref{eg:tabloid}.
A Hecke analogue of the theorem also holds: see~\S\ref{ss:ipf} below.

As pointed out by the referee, the recent paper~\cite{dn}
is concerned with constructing a basis for the annihilator 
of $\fieldc\tabloidl$ 
(and of its Hecke analogue $\perml$ whose definition is recalled below 
in~\S\ref{ss:permspecht}).    The answers are 
in terms of ``Murphy basis'',  which like the Kazhdan-Lusztig basis are 
known to be ``cellular''.   
\mysubsection{A question regarding the irreducible representations of the symmetric group}\mylabel{ss:iirrep}  
The question raised just above (in~\S\ref{ss:itabloid}) can be modified
to get one of more intrinsic appeal.  Given a partition~$\lambda$ of~$n$,
consider,  instead of the action of~$\symm_n$
on tabloids of shape~$\lambda$,   
the right cell 
module~$\rcmodl_\fieldc$ in the sense of Kazhdan-Lusztig~(\S\ref{ss:cellmod}),
or, equivalently (see~\S\ref{ss:mpiso}), the Specht 
module~$\spechtl_\fieldc$ (\S\ref{ss:specht}).   
The right cell modules are irreducible and every irreducible 
$\fieldc\symm_n$-module is isomorphic to $\rcmod(\mu)_\fieldc$ for 
some $\mu\partition n$
(\S\ref{ss:wstcon}).

The irreducibility of~$\rcmodl_\fieldc$ implies, by a well-known result of
Burnside (see, e.g., \cite[Chapter 8, \S4, No.~3, Corollaire~1]{bour}),
that the defining 
$\fieldc$-algebra map $\fieldc\symm_n\to\Endo_\fieldc{\rcmodlc}$ 
is surjective.    
The dimension of~$\rcmodlc$ (equivalently of~$\spechtl_\fieldc$)
equals the number~$\dl$
of standard tableaux of shape~$\lambda$ 
(\S\ref{sss:dlambda}, \S\ref{ss:wstcon}).  
Thus there exist $\dl^2$ elements of~$\fieldc\symm_n$,
even of~$\symm_n$ itself,  whose images in $\Endo_\fieldc{\rcmodlc}$
form a basis (for $\Endo_\fieldc{\rcmodlc}$).    We ask:
\begin{quote}
Is there is a natural
choice of such elements of~$\fieldc\symm_n$,  even of~$\symm_n$?
\end{quote}

Indeed there is,  as the following theorem says.  
As pointed out to us by Andrew Mathas, 
the theorem is a consequence 
of the {\em cellularity\/} in the sense of~\cite{glinv}
of the Kazhdan-Lusztig basis.  
\begin{theorem}[Consequence of cellularity 
of the Kazhdan-Lusztig basis~\cite{glinv}]
  \mylabel{t:irrep}  
Consider the Kazhdan-Lusztig basis elements of the group ring~$\fieldc\symm_n$ 
indexed by permutations of \shape~$\lambda$.   Their images
under the defining $\fieldc$-algebra map
$\fieldc\symm_n\to\Endo_\fieldc{\rcmodlc}$ form a basis 
for~$\Endo_\fieldc{\rcmodlc}$.
\end{theorem}
\noindent
By the Kazhdan-Lusztig basis elements of the group ring~$\fieldc\symm_n$,  
we mean the images in~$\fieldc\symm_n$ of
the Kazhdan-Lusztig basis elements of the Hecke algebra of~$\symm_n$
under the natural map setting the parameter value to~$1$ (\S\ref{s:setup}).
The \shape~of a permutation is defined using the 
RSK-correspondence~(\S\ref{ss:rsk}). 

The theorem is proved in~\S\ref{s:proof}.
   Some comments on the proof
of the theorem can be found in~\S\ref{ss:ipf}.

We do not know a natural choice of elements of the group~$\symm_n$ itself whose
images in~$\Endo_\fieldc{\rcmodlc}$ are a basis.    Permutations of~\shape\
$\lambda$ of course suggest themselves,   but they do
not in general have the desired property (Example~\ref{x:counter}).
\mysubsection{Comments on the proofs of 
Theorems~\ref{t:invariant}--\ref{t:irrep}}%
\mylabel{ss:ipf}
Properties of the Kazhdan-Lusztig basis of the Hecke algebra associated
to the symmetric group are the key to the proofs,
although the statements of Theorems~\ref{t:invariant} and~\ref{t:tabloid}
do not involve the Hecke algebra at all.   The relevant properties
are recalled in two instalments:  the first, in~\S\ref{s:setup},
is the more substantial;   the second, in~\S\ref{ss:pftinv}, consists 
of further facts
needed more specifically for the proof of Theorem~\ref{t:invariant}.

Theorem~\ref{t:irrep} follows by combining the
Wedderburn structure theory of semisimple algebras, as recalled
in~\S\ref{ss:wburn}, with the
following observation implicit in~\cite{geck} and explicitly formulated
in~\S\ref{ss:projections}:
\begin{quote}
  A Kazhdan-Lusztig $C$-basis element~$C_w$ 
kills the right cell (or equivalently Specht)
module corresponding to a shape~$\lambda$ unless $\lambda$ is
dominated by the \shape\ of the indexing permutation~$w$.
\end{quote}
The observation in turn follows easily from 
the combinatorial characterizations of the left, right, and
two-sided Kazhdan-Lusztig cells in terms of the RSK-corres\-pondence
and the dominance order on partitions (\S\ref{sss:cells}, \S\ref{sss:kllr}).
Our primary source for these characterizations, which are crucial to our
purpose,  is~\cite{geck}.

As pointed out to the authors by the referee among others, 
special properties of the Kazhdan-Lusztig basis of the Hecke algebra
of the symmetric group, as the one in the observation above,
have been noted and well studied.  In fact, they
have been axiomatized in~\cite{glinv},
where any basis enjoying these properties is termed {\em cellular\/}---see
also~\cite[Chapter~2]{mathasbook}.
Sections~\ref{s:proof} and~\ref{ss:encode} below are in effect working out
some consequences of cellularity;  and \ref{ss:cprod} is
in effect proving the cellularity of the Kazhdan-Lusztig basis using
results of~\cite{geck}.  Note that establishing cellularity is difficult,
there being a reliance on~\cite{geck} in our case.

Theorem~\ref{t:tabloid} follows by combining the above observation
with two well known facts: 
the isomorphism of the right cell module with the Specht module
and the well-known decomposition into irreducibles
of~$\fieldc\tabloid_\lambda$.
A Hecke analogue of Theorem~\ref{t:tabloid} also holds: see
Theorem~7 in the earlier version~\cite{arxiv:krs} of the present paper.  
A proof of it parallel to the proof of Theorem~\ref{t:tabloid} as in
here can be given using results of~\cite{djqschur}.    The proof
in~\cite{arxiv:krs} is different and more in keeping
with the ideas developed here.

The main technical point in the proof of Theorem~\ref{t:invariant} is 
isolated as Lemma~\ref{l:cfreemlin},  which is a two sided analogue
of~\cite[Lemma~2.11]{mp} recalled below as Proposition~\ref{p:pftinv2}.
\mysubsection{Analogue of Theorem~\ref{t:irrep} over fields of positive
characteristic}\mylabel{ss:iposchar}
The Hecke algebra and its Kazhdan-Lusztig basis make sense over
an arbitrary base~(\S\ref{s:setup}).    The cell modules and
Specht modules are also defined and 
isomorphic over any base~(\S\ref{ss:cellmod},
\S\ref{ss:specht}, \S\ref{ss:mpiso}).     
Thus we can ask for the analogue of
Theorem~\ref{t:irrep} over a field of arbitrary characteristic, keeping
in mind of course that the cell modules may not be irreducible any longer.
We prove:
\begin{theorem}
  \mylabel{t:iposchar}
Let~$\field$ be a field of positive characteristic~$p$.   Let~$\lambda$
be a partition of a positive integer~$n$ no part of which is repeated
$p$~or more times.   Suppose that the right cell module~$\rcmodl_\base$ is
irreducible.   
Consider the Kazhdan-Lusztig basis elements of the group ring~$\field\symm_n$ 
indexed by permutations of \shape~$\lambda$.   Their images
under the defining $\field$-algebra map
$\field\symm_n\to\Endo_\field{\rcmodl_\base}$ form a basis 
for~$\Endo_\field{\rcmodl_\base}$.
\end{theorem}
\noindent
The theorem is a special case of Theorem~\ref{t:cpklbasis} proved
in~\S\ref{s:poschar}.    Like the proof of Theorem~\ref{t:irrep},
that of Theorem~\ref{t:cpklbasis} too uses the observation
formulated in~\S\ref{ss:projections},  but, the group ring~$\base\symm_n$
being not necessarily semisimple,  we cannot rely on Wedderburn structure
theory any more.   Instead we take a more head-on approach:
%
\begin{quote}%
Choosing a convenient basis of~$\Endo_\base\rcmodl_\base$,
we express as linear combinations of these basis elements the images 
in~$\Endo_\base\rcmodl_\base$ of the appropriate
Kazhdan-Lusztig basis elements of~$\base\symm_n$.  Denoting
by~$\bigmatrixl$ the resulting square matrix 
of coefficients,  we give an explicit formula for its determinant
$\det\bigmatrixlk$.
\end{quote}
In fact,  we obtain a formula for~$\det\bigmatrixl$,  where 
$\bigmatrixl$ is the analogous matrix of coefficients over an arbitrary base and
over the Hecke algebra (rather than the group ring): see~\S\ref{s:matrixgl} for
details.
We then need only specialize to get $\det\bigmatrixl_\base$.
Given the formula, it is a relatively easy matter to get a criterion
for $\det\bigmatrixl_\base$ not to vanish,  thereby proving
Theorem~\ref{t:cpklbasis}.   
\mysubsection{A hook length formula for the determinant of~$\matrixgl$}%
\mylabel{ss:icomputation}
To obtain the formula for $\det\bigmatrixl$,  
we first show that~$\bigmatrixl$
has a nice form which enables us to reduce the computation to that
of the determinant of a matrix~$\matrixgl$ of much smaller size.
We discuss how this is done.

The basis of~$\Endo\rcmodl$ with respect to which the matrix~$\bigmatrixl$
is computed suggests itself:  $\rcmodl$ has a 
basis consisting of classes of Kazhdan-Lusztig elements $C_w$,
where $w$ belongs to a right cell of shape~$\lambda$ 
of~$\symm_n$ (\S\ref{ss:cellmod});  considering the endomorphisms
which map one of these basis elements to another (possibly the same) 
and kill the rest,
we get the appropriate basis for $\Endo\rcmodl$.   This means that
the matrix~$\bigmatrixl$ encodes the multiplication table for
Kazhdan-Lusztig basis elements~$C_w$ indexed by permutations of
\shape~$\lambda$, modulo those indexed by permutations of lesser
shape in the dominance order.   

The special (cellularity) properties of the
Kazhdan-Lusztig elements now imply that the matrix~$\bigmatrixl$,
which is of size $\dl^2\times\dl^2$ (where $\dl$ is the number of
standard tableaux of shape~$\lambda$),
is a `block scalar' matrix,  i.e.,  when broken up into blocks
of size $\dl\times\dl$,  only the diagonal blocks are non-zero,
and all the diagonal blocks are equal.  Denoting by~$\matrixgl$ the
diagonal block,  we are thus reduced to computing the determinant
of~$\matrixgl$.    The details of this reduction are worked out
in~\S\ref{s:matrixgl}.   

The formula for the determinant 
of~$\matrixgl$ is given in Theorem~\ref{t:formula}, the main
ingredients in the proof of which are formulas from~\cite{djblocks} 
and~\cite{jmathas}.  The relevance
of those formulas to the present context is not clear at first sight.
They are about the determinant, denoted~$\gramdet$, of the matrix of a certain
bilinear form,  the Dipper-James form,  on the Specht module~$\spechtl$, 
computed with respect to the `standard basis' of~$\spechtl$;  while
$\matrixgl$ has to do with multiplication of Kazhdan-Lusztig basis elements.
The connection between~$\gramdet$ and $\det\matrixgl$ is established
in~\S\ref{s:gandgram} (see Equation~(\ref{e:gandgram})) using 
results of~\cite{mp}.
\mysubsection{On the irreducibility of Specht modules}%
\mylabel{ss:icarter}   Finally, we discuss another application of
the formula for the determinant of the 
matrix~$\matrixgl$ introduced in~\S\ref{ss:icomputation}.
Suppose that the determinant did not vanish when the Hecke algebra
is specialized to group ring and the scalars extended to a
field~$\base$.    Then, evidently, the images in~$\Endo_\base\rcmodl_\base$ of
the Kazhdan-Lusztig basis elements~$C_w$,
as $w$ varies over permutations of \shape~$\lambda$,  form a basis
for $\Endo_\base\rcmodl_\base$,  which means in particular that the
defining map $\base\symm_n\to\Endo_\base\rcmodl_\base$ is surjective, and
so~$\rcmodl_\base$ is irreducible.    

In other words,  the non-vanishing
of~$\det\matrixgl$ in~$\base$ gives a criterion for the irreducibility
of~$\rcmodl_\base$ (equivalently, of~$\spechtl_\base$).    The criterion
thus obtained matches precisely the one conjectured by Carter and proved
in~\cite{james:carter,jmathas}.   
We thus obtain an independent proof of the Carter criterion.
The details are worked out in~\S\ref{s:poschar}.
\mysubsection{Acknowledgments}  Thanks to the {\sffamily GAP\/} 
program, computations performed on which were useful in 
formulating the results;  to the 
{\sffamily Abdus Salam International Centre
for Theoretical Physics\/},  during a visit to which of one
of the authors much of this work was brought to completion;
to the {\sffamily Skype\/} program, using which the authors
were able to stay in touch during that visit.

Thanks especially to
John Graham, Andrew Mathas, and the referee for their comments.

%% file: comb.tex
\mysection{Recall of some basic notions}
\mylabel{s:comb}
\noindent
We recall in this section the basic combinatorial and representation
theoretic notions that we need.   Note that our definition 
of the RSK-correspondence (\S\ref{ss:rsk}) differs from the standard
(as e.g. in~\cite[Chapter~4]{fulton}) by a flip.    

Throughout~$n$ denotes a positive integer.
\mysubsection{Partitions and shapes}\mylabel{ss:partshape}
By a {\em partition\/}~$\lambda$ of~$n$, written $\lambda\partition n$,
is meant a sequence $\lambda_1\geq
\ldots \geq \lambda_r$ of positive integers such
that $\lambda_1+\ldots+\lambda_r=n$.   The integer~$r$ is 
the {\em number of parts\/} in~$\lambda$.  We often write
$\lambda=(\lambda_1,\ldots,\lambda_r)$;  sometimes even 
$\lambda=(\lambda_1,\lambda_2,\ldots)$.   When the latter notation
is used,  it is to be understood that $\lambda_t=0$ for $t>r$.

Partitions of~$n$ are
in bijection with {\em shapes of Young diagrams\/} (or simply {\em shapes\/}) 
with $n$ boxes:
the partition $\lambda_1\geq\ldots\geq\lambda_r$ corresponds to the
shape with $\lambda_1$ boxes in the first row, $\lambda_2$~in the
second row, and so on,  the boxes being arranged left- and top-justified.
Here for example is the shape corresponding to the partition~$(4,3,1)$ of~$8$:
\[\begin{array}{|c|c|c|c|}
\hline
& & & \\ 
\hline
& & \\
\cline{1-3} \\
\cline{1-1}
\end{array}\]
Partitions are thus identified with shapes and the two terms are used
interchangeably.
\mysubsection{Dominance order on partitions}\mylabel{ss:domination}
Given partitions $\mu =
(\mu_1, \mu_2, \ldots)$ and $\lambda = (\lambda_1, \lambda_2,
\ldots)$ of~$n$, we say $\mu$ {\em dominates\/}
$\lambda$, and write $\mu \domgeq \lambda$, if \[ 
\mu_1 \geq \lambda_1,\quad\mu_1+
\mu_2 \geq \lambda_1 + \lambda_2,\quad 
\mu_1 + \mu_2 +\mu_3 \geq \lambda_1 + \lambda_2+\lambda_3,\quad\ldots.\] 
We write $\mu\domr\lambda$ if $\mu\domgeq\lambda$ and $\mu\neq\lambda$.
The partial order $\domgeq$ on the set of
partitions (or shapes) of $n$ will be referred to as the {\em
  dominance order\/}.
\mysubsubsection{The partition~$\lnd$}\mylabel{sss:lnd}
Given integers~$n$ and~$d$,   there exists a unique partition~$\lnd\partition n$
that has at most~$d$ parts and is smallest in the dominance order among
those with at most~$d$ parts.    For example, $\lambda(8,3)=(3,3,2)$.
\mysubsection{Tableaux and standard tableaux}\mylabel{ss:youngtab}
A {\em Young tableau\/}, or just {\em tableau\/}, of shape~$\lambda\partition n$
is an arrangement of the numbers $1$, \ldots, $n$  in the boxes 
of shape~$\lambda$.   There are, evidently, $n!$ tableaux of shape~$\lambda$.
A tableau is {\em row standard\/} (respectively,
{\em column standard\/}) if in every row (respectively, column)
the entries are increasing left to right (respectively, top to bottom).
A tableau is {\em standard\/} if it is both row standard and column
standard.
An example of a standard tableau of shape $(3, 3, 2)$:
\[\begin{array}{|c|c|c|}
\hline
1&3 &5 \\ 
\hline 
2&6 &8 \\
\hline 4 &7 \\
\cline{1-2}
\end{array}\]
\mysubsubsection{The number of standard tableaux}\mylabel{sss:dlambda}
The number of standard tableaux of a given shape~$\lambda\partition n$ 
is denoted~$\dl$.
There is a well-known `hook length formula' for it~\cite{frt}:   
$\dl=n!/\prod_\beta h_\beta$,  where $\beta$ runs over all boxes of 
shape~$\lambda$ and
$h_\beta$ is the {\em hook length\/} of the box~$\beta$ which is defined
as one more than the sum of the number of boxes to the right of~$\beta$
and the number of boxes below~$\beta$.   

The hook lengths for shapes $(3,3,2)$ and $(4,3,1)$ are shown below:
\[\begin{array}{|c|c|c|}
\hline
5 &4 &2 \\ 
\hline 
4&3 &1 \\
\hline 2&1 \\
\cline{1-2}
\end{array}\quad\quad\quad\quad\quad\quad
\begin{array}{|c|c|c|c|}
\hline
6 &4 &3 & 1\\ 
\hline
4 &2 & 1\\
\cline{1-3} 1\\
\cline{1-1}
\end{array}
\]
Thus $d(3,3,2)=8!/(5\cdot4\cdot2\cdot4\cdot3\cdot1\cdot2\cdot1)=42$ and $d(4,3,1)=8!/6\cdot4\cdot3\cdot1\cdot4\cdot2\cdot1\cdot1=70$.
\mysubsection{The RSK-correspondence and the \shape\ of a permutation}%
\mylabel{ss:rsk} 
The {\em Robinson-Schensted-Knuth correspondence\/} ({\em RSK
correspondence\/} for short) is
a well-known procedure that sets up a bijection between 
the symmetric group $\symmn$ and ordered pairs of standard tableaux 
of the same shape with $n$ boxes.   
%
We do not recall here the
procedure, referring the reader instead to~\cite[Chapter~4]{fulton}.
It will be convenient for our purposes to modify slightly the procedure
described in~\cite{fulton}.  \par
Denoting by $(A(w),B(w))\leftrightarrow w$
the bijection of~\cite{fulton},    what we mean by {\em 
RSK correspondence\/} is the bijection $(B(w),A(w))\leftrightarrow 
w$;   since $A(w)=B(w^{-1})$ and $A(w^{-1})=B(w)$ (see \cite[Corollary on page~41]{fulton}),
we could equally well define our RSK correspondence as
$(A(w),B(w))\leftrightarrow w^{-1}$.
The {\em \shape\/} of a permutation~$w$ is defined to be the shape of
either of~$A(w)$, $B(w)$.
\mysubsubsection{An example}\mylabel{sss:xrsk}   The permutation
$(1542)(36)$ (written as a product of disjoint cycles) has \shape\
$(3,2,1)$.     Indeed it is mapped under the RSK correspondence in our sense
to the ordered pair $(A,B)$ of standard tableaux,  where:
\[A= 
\begin{array}{|c|c|c|}
\hline
1&3 &5 \\ 
\hline 
2&4 \\
\cline{1-2} 6 \\
\cline{1-1}
\end{array}
\quad\quad\quad
B=
\begin{array}{|c|c|c|}
\hline
1&2 &3 \\ 
\hline 
4&6 \\
\cline{1-2} 5 \\
\cline{1-1}
\end{array}
\]
\mysubsubsection{Remark}\mylabel{sss:rskremark}
The justification for our modification of the standard definition of 
RSK-correspondence is that it was the simplest way we could
think up of reconciling the notational conflict among the two sets of papers
upon which we rely:~\cite{dj,mp} and~\cite{geck}.
Permutations act on the 
right in the former---a
convention which we too follow---but 
on the left in the latter.     The direction in which they act makes
a difference to statements involving the RSK correspondence: most
importantly for us, to the characterization of one-sided cells
(see~\S\ref{sss:cells} below).     This creates a problem:  we cannot
be quoting literally from both sets of sources without changing something.
Altering the definition of RSK correspondence as above 
is the path of least resistance,  and allows us to quote more or less
verbatim from both sets.

\mysubsection{Tabloids and tabloid representations}\mylabel{ss:tabloid}
Let $\lambda=(\lambda_1,\lambda_2,\ldots)\partition n$.
A {\em tabloid\/} of shape~$\lambda$ 
is a partition of the set~$[n]:=\{1,\ldots,n\}$ into an ordered $r$-tuple 
of subsets, the first consisting of $\lambda_1$ elements, the second
of $\lambda_2$ elements, and so on.    Depicted below are two tabloids
of shape $(3,3,2)$:
\[\begin{array}{ccc}
\hline
1 & 3 & 5\\ 
\hline 
7 & 8 & 9 \\
\hline 4 & 6 \\
\cline{1-2}
\end{array} \quad\quad\quad\quad\quad\quad
\begin{array}{ccc}
\hline
3 & 5  & 8 \\ 
\hline 
1 & 6 & \\
\hline 2 & 7 \\
\cline{1-2}
\end{array}\]
The members of the first subset are arranged in increasing order in the
first row,  those of the second subset in the second row, and so on.

Given a tableau~$T$ of shape~$\lambda$,  it determines, in the obvious
way,  a tabloid of shape~$\lambda$ denoted $\{T\}$:   
the first subset consists of the elements in the first row,
the second of those in the second row, and so on.

The defining action
of~$\symm_n$ on $[n]$ induces, in the obvious way,  an action on the
set $\tabloidl$ of tabloids of shape~$\lambda$.
The free $\integers$-module~$\itabl$ with $\tabloidl$
as a $\integers$-basis provides therefore a 
linear representation of~$\symm_n$ over~$\integers$.
By base change we get such a representation over
any commutative ring with 
unity~$\field$: $\field\tabloidl:=\integers\tabloidl\tensor_\integers\field$.  
We call it the {\em tabloid representation\/} corresponding to 
the shape~$\lambda$.
\mysubsection{Specht modules}\mylabel{ss:specht}
\newcommand\scripte{\mathfrak{e}}
The Specht module corresponding to a partition~$\lambda\partition n$ 
is a certain $\symm_n$-submodule of the tabloid 
representation~$\integers\tabloidl$ just defined.
For a tableau~$T$ of shape~$\lambda$, 
define $\scripte_T$ in $\integers\tabloidl$ by
\[ \scripte_T := \sum \sign(\sigma)\{T\sigma\} \] 
where the sum is taken over permutations $\sigma$ of~$\symm_n$ 
in the column stabiliser of~$T$, $\sign(\sigma)$ denotes the
sign of~$\sigma$, and $\{T\sigma\}$ denotes the tabloid corresponding
to the tableau~$T\sigma$ in the obvious way (see~\S\ref{ss:tabloid}).
The {\em Specht module $S^{\lambda}$\/} is the linear
span of the $\scripte_T$ as $T$ runs over all
tableaux of shape $\lambda$.   It is an~$\symm_n$-submodule 
of~$\integers\tabloidl$ with $\integers$-basis
$\scripte_T$, as $T$ varies over standard tableaux 
(see, for example,~\cite[\S7.2]{fulton}). 
By base change we get the Specht 
module~$\spechtl_\base$ over any commutative ring with identity~$\base$:
$\spechtl_\base:=\spechtl\tensor_\integers\base$.
Evidently, $\spechtl_\base$ is a free $\base$-module of rank the number~$\dl$
of standard tableaux of shape~$\lambda$ (\S\ref{sss:dlambda}).

%% file: hecke.tex
\mysection{Set up:  Hecke algebra and Kazhdan-Lusztig cells}\mylabel{s:setup}
\noindent
Let~$n$ denote a fixed positive integer and~$\symm_n$ the symmetric
group on~$n$ letters.
Let $S$ be the subset consisting of the simple transpositions
$(1,2)$, $(2,3)$, \ldots, $(n-1,n)$
of the symmetric group $\symm_n$.   Then $(\symm_n, S)$ is a Coxeter
system in the sense of~\cite[Chapter~4]{bour}.
Let $A:=\integers[v,v^{-1}]$, the Laurent polynomial ring in
the variable~$v$ over the integers.   
\mysubsection{The Hecke algebra and its $T$-basis}\mylabel{ss:hecke}
Let $\hecke$ be the {\em Iwahori-Hecke algebra\/}
corresponding to~$(\symm_n, S)$, with notation as in~\cite{geck}.
Recall that $\hecke$ is an $A$-algebra:
it is a free $A$-module with basis $T_w$, $w\in\symm_n$, the
multiplication being defined by 
\[
T_sT_w=\left\{
\begin{array}{ll}
T_{sw} & \textup{if $\length(sw)=\length(w)+1$}\\
(v-v^{-1})T_w+T_{sw} & \textup{if $\length(sw)=\length(w)-1$}\\
\end{array}
\right.
\]
for $s\in S$ and $w\in\symmn$ and $\ell$ is the length function.
We put \[\epsilon(w):=(-1)^{\length(w)}\quad\textrm{and}
\quad v_w:=v^{\length(w)}\quad \textrm{for $w\in\symm_n$}.\]

An induction on length gives (in any case, see 
\cite[Lemma~2.1~(iii)]{dj} for (\ref{eq:tutu'})):
\begin{gather}
  \label{eq:twtw'}
T_wT_{w'}=T_{ww'}\quad\quad\quad\textup{if $\length(w)+\length(w')=\length(ww')$.}\\
  \label{eq:tutu'}
T_u T_{u'}= T_{uu'} + \sum_{uu'<w}
a_w T_w \quad\quad\textup{for $u$, $u'$ in~$\symm_n$}
  \end{gather}   
Here, as elsewhere, $<$ denotes the Bruhat-Chevalley partial order on~$\symm_n$.
In particular,  the coefficient of~$T_1$ in $T_u T_{u'}$ is
non-zero if and only if $u'=u^{-1}$ and equals~$1$ in that case.
\mysubsubsection{The relation between~$v$ and~$q$}\mylabel{sss:vandq}
We follow the conventions
of~\cite{geck}.    In particular,   to pass from our notation to that
of~\cite{kl}, \cite{dj}, or \cite{mp},   
we need to replace $v$ by $q^{1/2}$ and $T_w$ by 
$q^{-\length(w)/2}T_w$.    
\mysubsubsection{Specializations of the Hecke algebra}%
\mylabel{sss:hspecial}
Let $\base$ be a commutative ring with unity and~$a$ an invertible
element in~$\base$.
There is a unique ring
homomorphism~$A\to\base$ defined by~$v\mapsto a$.
We denote by~$\hecke_\base$ the $\base$-algebra~$\hecke\tensor_A\base$
obtained by extending the scalars to~$\base$ via this homomorphism.  
We have a natural $A$-algebra homomorphism $\hecke\to\hecke_\base$ given
by $h\mapsto h\tensor1$.
By abuse of notation,  we continue to use the same symbols for the
images in $\hecke_\base$ of elements of~$\hecke$ as for those elements
themselves. 
If~$M$ is a (right)
$\hecke$-module,  $M\tensor_A\base$ is naturally a (right) 
$\hecke_\base$-module.   

An important special case is when we take $a$ to be the unit element~$1$
of~$\base$.
We then have a natural identification of $\hecke_\base$ with
the group ring $\base\symm_n$, under which $T_w$ maps
to the permutation~$w$ in~$\base\symm_n$.    

Regarding the semisimplicity 
of~$\hecke_\base$, we have this result~\cite[Theorem~4.3]{djblocks}:
\begin{quote}
Assuming $\base$ to be a field,
$\hecke_\base$ is semi-simple except precisely when
\begin{itemize}
\item 
either $a^2=1$ and the characteristic of~$\base$ is $\leq n$
\item
or $a^2\neq1$ is a primitive $r^{\textup{th}}$ root of unity for some
$2\leq r\leq n$.   
\end{itemize}
\end{quote}

\mysubsubsection{Two ring involutions and an $A$-antiautomorphism}%
\mylabel{sss:2inv1aa}
We use the following two involutions on~$\hecke$ both of which extend
the ring involution $a\mapsto\bar{a}$ of~$A$ defined by 
$v\mapsto\bar{v}:=v^{-1}$:  
\[ 
\bar{\sum a_w T_w}:=\sum \bar{a_w} T^{-1}_{w^{-1}}
\quad\quad\quad
j(\sum a_wT_w):= 
\sum \epsilon_w \bar{a_w}T_w \quad\quad\quad
\textup{($a_w\in A$)}
\]
These commute with each other and so their composition,
denoted $h\mapsto h^\dag$,  is an $A$-algebra involution of~$\hecke$.
The $A$-algebra anti-automorphism~$(\sum a_w T_w)^*=
\sum a_w T_{w^{-1}}$ allows passing back and forth
between statements about left cells and orders and those
about right ones (\S\ref{sss:cells}).
\mysubsection{Kazhdan-Lusztig $C'$- and $C$-basis}\mylabel{ss:klbasis}
Two types of $A$-bases for~$\hecke$ are introduced 
in~\cite{kl},  denoted $\{C'_w\st w\in \symmn\}$ and $\{C_w\st w\in \symmn\}$.
They
are uniquely determined by the respective 
conditions~\cite[Theorem~5.2]{lusztig}:
\begin{equation}\label{eq:cc'char}\begin{split}
\bar{C'_w} & = C'_w \quad \textup{ and }\quad C_w'\equiv T_w\mod{\hecke_{<0}}\\
\bar{C_w}& = C_w \quad \textup{ and }\quad
C_w\equiv T_w \mod{\hecke_{>0}}
\end{split}\end{equation}
where
\[\textup{
$\hecke_{<0}:=\sum_{w\in \symmn} A_{<0}T_w$, $A_{<0}:=v^{-1}\integers[v^{-1}]$}
\quad\quad
\textup{
$\hecke_{>0}:=\sum_{w\in \symmn}
A_{>0}T_w$, $A_{>0}:=v\integers[v]$}\]
The anti-automorphism $h\mapsto h^*$ and the ring involution
$h\mapsto \bar{h}$ commute 
with each other, so that, by the characterization~(\ref{eq:cc'char}):
\begin{equation}
  \label{eq:c'star}
  (C_x)^*=C_{x^{-1}}\quad\quad\quad\quad\quad\quad (C_x')^*=C_{x^{-1}}'
\end{equation}
We have by~\cite[Theorem~1.1]{kl}:
\begin{equation}\label{e:ct}
C'_w=T_w+\sum_{y\in \symmn, y<w} p_{y,w}T_y
\quad\quad\quad
C_w=T_w+\sum_{y\in \symmn, y<w} \epsilon_y\epsilon_w \bar{p_{y,w}}T_y
\end{equation}
where 
$<$ denotes the Bruhat-Chevalley order on~$\symmn$, and $p_{y,w}\in A_{<0}$
for all $y<w$,  from which it is clear that 
\begin{equation}
  \label{eq:cjc'}
C_w=\epsilon_w j(C_w')
\end{equation}
Combining (\ref{eq:cc'char}) with (\ref{eq:cjc'}),  we obtain
\begin{equation}
  \label{eq:c'dag}
  C_w=\epsilon_w (C'_w)^\dag
\end{equation}
\mysubsubsection{Notation}\mylabel{sss:notation}
For a subset $\subsets$ of~$\symm_n$,   denote by
$\langle C_y\st y\in\subsets\rangle_A$ the $A$-span in~$\hecke$ of 
$\{C_y\st y\in\subsets\}$.   
For an $A$-algebra~$\base$,
denote by $\langle C_y\st y\in\subsets\rangle_\base$ 
the $\base$-span in~$\hecke_\base$ of
$\{C_y\st y\in\subsets\}$.
 Similar meanings are attached to
$\langle T_y\st y\in\subsets\rangle_A$ and
$\langle T_y\st y\in\subsets\rangle_\base$.
\mysubsubsection{A simple observation}\mylabel{sss:observe}
From~(\ref{e:ct}), we get 
$T_w\equiv C_w\bmod \langle T_x\st x<w\rangle_A$.   From this in turn
we get, by induction on the Bruhat-Chevalley order, the following:  
for a subset~$\subsets$ 
of~$\symm_n$,   
the (images of) elements~$T_w$, $w\in\symm_n\setminus\subsets$,
form a basis for the $A$-module $\hecke/\langle
C_x\st x\in\subsets\rangle_A$.    The same thing holds also in
specializations $\hecke_\field$ of~$\hecke$ (\S\ref{sss:hspecial}):
the (images of) elements~$T_w$, $w\in\symm_n\setminus\subsets$,
form a basis for the $\field$-module
 $\hecke_\field/\langle
C_x\st x\in\subsets\rangle_\field$.
%
\mysubsection{Kazhdan-Lusztig orders and cells}\mylabel{ss:klorders}
Let $y$ and $w$ in $\symm_n$. Write $y\klla w$ if,
for some element~$s$ in~$S$, the coefficient of $C_y$ is non-zero 
in the expression of $C_s C_w$ as a $A$-linear combination of 
the basis elements~$C_x$.
Replacing all occurrences of `$C$' by `$C'$' in this definition
would make no difference.
The Kazhdan-Lusztig {\em left pre-order\/}
is defined by: $y\kll w$ if there exists 
a chain $y=y_0\klla\cdots\klla y_k=w$;   the {\em left equivalence\/} relation
by: $y\kllsim w$ if $y\kll w$ and $w\kll y$.
Left equivalence classes are called {\em left cells\/}.
Note that $\sum_{x\kll w}A C_x$ is a left ideal containing the left
ideal~$\hecke C_w$.

Right pre-order, equivalence, and cells are defined similarly.
The {\em two sided pre-order\/} is defined by:  $y\kllr w$ if there
exists a chain $y=y_0$, \ldots, $y_k=w$  such that,  for $0\leq j<k$,
either $y_j\kll y_{j+1}$ or $y_j\klr y_{j+1}$.  Two sided equivalence
classes are called {\em two sided cells.}
\mysubsection{Cells and RSK Correspondence}\mylabel{ss:cellsandkrs}
We now recall the combinatorial characterizations of one and two sided cells in
terms of the RSK correspondence (\S\ref{ss:rsk}) and the dominance
order on partitions (\S\ref{ss:domination}).   These statements
are the foundation on which this paper rests.   The ones
in~\S\ref{sss:kllr},~\ref{sss:kllr} are used repeatedly,  but
the more subtle one in~\S\ref{sss:unrelate} is used only once, namely 
in the proof of Theorem~\ref{t:invariant}:  it is used in the proof of
Lemma~\ref{l:cfreemlin} which is the main ingredient
in the proof of that theorem.

Write $(P(w),Q(w))$ for the ordered pair of standard Young tableaux 
associated to a permutation~$w$ by the RSK correspondence (in our sense---see
\S\ref{ss:rsk}).   
Call $P(w)$ the
{\em $P$-symbol\/} and $Q(w)$ the {\em $Q$-symbol\/} of~$w$.
It will be convenient to use such notation as 
$(P(w),Q(w))$ for the permutation~$w$, 
$C_{(P(w),Q(w))}$ or $C(P(w),Q(w))$ for 
the Kazhdan-Lusztig $C$-basis element~$C_w$.
\mysubsubsection{Cells in terms of symbols}\mylabel{sss:cells}
Two permutations  
are left equivalent if and only if they have the
same $Q$-symbol;  right equivalent if and only if 
the same $P$-symbol;  two sided equivalent if and only if 
the same \shape.
See~\cite[Corollary~5.6]{geck} (and comments therein about \cite[\S5]{kl}, 
\cite{ariki}).
\mysubsubsection{The $\kllr$ relation in terms of dominance}%
\mylabel{sss:kllr}
We have 
$y\kllr w$ if and only if
$\krsshape{(y)}\domleq \krsshape{(w)}$,
where $\domleq$ is the usual dominance order on partitions:  
$\lambda\domleq\mu$ if
$\lambda_1\leq \mu_1$, $\lambda_1+\lambda_2
\leq \mu_1+\mu_2$, \ldots.
See~\cite[Theorem~5.1]{geck} (and comments therein about~\cite[2.13.1]{dps}).
We write $\lambda\doml\mu$ for $\lambda\domleq\mu$ and $\lambda\neq\mu$.
\mysubsubsection{Unrelatedness of distinct one sided cells 
in the same two sided cell}%
\mylabel{sss:unrelate}   
If $x\kll y$ and $x\kllrsim y$,  then $x\kllsim y$.
See~\cite[Theorem~5.3]{geck} (and comments therein 
about~\cite[Lemma~4.1]{lu81}).
\mysubsection{Cell modules}\mylabel{ss:cellmod}
It follows from the definition of the pre-order $\kll$ that
the $A$-span $\langle C_y\st y\kll w\rangle_A$ of
$\{C_y\st y\kll w\}$,  for $w$ in~$\symm_n$ fixed,
is a left ideal of~$\hecke$;  so is $\langle C_y\st y\kllt w\rangle_A$.
The quotient $\lcmod(w):=\langle C_y\st y\kll w\rangle_A/
\langle C_y\st y\kllt w\rangle_A$ is called the {\em left cell
module\/} associated to~$w$.   It is a left $\hecke$-module.
Right cell modules $\rcmod(w)$ and two sided cell modules are defined
similarly.   They are right modules and bimodules respectively.

Let $y$ and $w$ be permutations of the same \shape~$\lambda$.
The left cell modules $L(y)$ and $L(w)$ are then $\hecke$-isomorphic.
In fact, the association
$C_{(P,Q(y))}\leftrightarrow C_{(P,Q(w))}$ gives an
isomorphism:  
see~\cite[\S5]{kl}, \cite[Corollary~5.8]{geck}.   
The right cell modules $\rcmod(y)$ and $\rcmod(w)$ are
similarly isomorphic, and we sometimes write $\rcmodl$ for $\rcmod(y)
\simeq\rcmod(w)$.

When a homomorphism from $A$ to a commutative ring~$\field$ is specified,
such notation as $\rcmod(w)_\field$ and $\rcmodl_\field$ make 
sense: see~\S\ref{sss:hspecial}.
\mysubsection{A key observation regarding 
images of $C$-basis elements in endomorphisms of cell modules}
\mylabel{ss:projections}
The image of~$C_y$ in $\Endo{\rcmodl}$ vanishes unless
$\lambda\domleq\textup{\shape}(y)$,  for, if 
$C_z$ occurs
with non-zero coefficient in~$C_xC_y$ (when expressed as an $A$-linear
combinations of the $C$-basis),  where $\textup{\shape}(x)=\lambda$,  and
$\lambda\not\domleq\textup{\shape}(y)$,
then $z\kll y$ (by definition),  so $\textup{\shape}(z)
\domleq y$ (\S\ref{sss:kllr}),  which means that $\textup{\shape}(z)
\neq\lambda$, so $z\not\sim_R x$~(\S\ref{sss:cells}).

%% file: invapp.tex
\mysection{Applications to invariant theory}\mylabel{s:invapp}
\noindent
In~\S\ref{ss:pftinv}, Theorem~\ref{t:invariant} stated in~\S\ref{ss:iinv}
is proved.  
In the later subsections,
applications of the theorem to rings of multilinear and 
polynomial invariants are discussed. The base
$\base$ is an arbitrary commutative ring with unity in~\S\ref{ss:pftinv}
but in the later subsections it is assumed to satisfy further hypothesis.
\mysubsection{
Proof of Theorem~\ref{t:invariant}}\mylabel{ss:pftinv}
Our goal in this subsection is to prove Theorem~\ref{t:invariant} stated
in~\S\ref{ss:iinv}.
Let $n$ and $d$ be positive integers, $d<n$.  (The theorem clearly
holds when~$d\geq n$.) The main ingredient of the proof is 
Lemma~\ref{l:cfreemlin} below.   Once the lemma is proved, the
theorem itself follows easily: see~\S\ref{sss:pftinv}.    For the lemma
we need some combinatorial preliminaries, beyond those recalled 
in~\S\ref{s:comb}.
\mysubsubsection{Preliminaries to the proof}%
\mylabel{sss:invprelims}
For $\lambda$ a partition of~$n$,
\begin{itemize}
\item
$\lambda'$ denotes the {\em transpose\/} of $\lambda$.  E.g.,
$\lambda'=(3,2,2,1)$ for $\lambda=(4,3,1)$.
\item $\tupl$ denotes the standard tableau of shape~$\lambda$ in which the
numbers $1$,~$2$,~\ldots,~$n$ appear in order along successive rows;
$\tlowl$ is defined similarly, with  
`columns' replacing `rows'.
E.g., for $\lambda=(4,3,1)$, we have:
\[
\tupl=\begin{array}{|c|c|c|c|}
\hline
1 & 2 & 3 & 4\\ 
\hline 
 5 & 6 &7 \\
\cline{1-3}
8 \\
\cline{1-1}
\end{array}\quad\quad
\tlowl=\begin{array}{|c|c|c|c|}
\hline
1 & 4 & 6 & 8\\ 
\hline 
 2 & 5 &7 \\
\cline{1-3}
3 \\
\cline{1-1}
\end{array}\quad\quad
\]
\item $\Wl$ denotes the row stabilizer of~$\tupl$.  It is a parabolic subgroup
of~$\symm_n$.    E.g., for $\lambda=(4,3,1)$,  $\Wl$ is isomorphic to
the product~$\symm_4\times\symm_3\times\symm_1$.
\item $\wnotl$ denotes the longest element of~$\Wl$.  E.g., when $n=8$ and
$\lambda=(4,3,1)$,  the sequence $(1\wnotl, \ldots, n\wnotl)$ is
$(4,3,2,1,7,6,5,8)$.
\item $\dija:=\{w\in\symm_n\st\textup{$t^\lambda w$ is row
     standard}\}$.  Clearly $\dija$ is a set of right coset
   representatives of $\Wl$ in $\symm_n$ and  $\wl$ is an element
   of~$\dija$.
\end{itemize}
\begin{proposition}
  \mylabel{p:pftinv} 
For $\lambda$ a partition of~$n$,
\begin{enumerate}
\item\label{i:lwd} $\length(wd)=\length(w)+\length(d)$,
for $w\in\Wl$ and $d\in\dija$.
\item\label{i:dindija} $d\in\dija$ is the unique element of minimal length in $\Wl d$.
\item\label{i:wnotldija}
$\wnotl\dija=\{w\st w\klr\wnotl\}$.
Thus $w\klr\wnotl$ if and
only if in every row of~$\tupl w$ the entries are decreasing to
the right. 
 \item\label{i:wnotl}
 $\wnotl$ is of shape~$\lambdap$: it corresponds under RSK
 to $(t_{\lambdap},t_{\lambdap})$.
\end{enumerate}
\end{proposition}
\begin{proof}
(\ref{i:lwd}) and (\ref{i:dindija}) are elementary to see: in any case,
see \cite[Lemma~1.1~(i),~(ii)]{dj}.
(\ref{i:wnotl})~is evident.  For~(\ref{i:wnotldija}) see~\cite[\S2.9]{mp}. 
\end{proof}

   Now, fix notation as in~\S\ref{s:setup}:
$A$ denotes the ring~$\integers[v,v^{-1}]$, $\hecke$ the Hecke algebra,
$C_w$ the Kazhdan-Lusztig basis element corresponding to the permutation~$w$,
etc.     
\begin{proposition}
  \mylabel{p:pftinv2} \textup{(\cite[Lemma~2.11]{mp})}
The $A$-span $\langle C_w\st w\klr\wnotl\rangle_A$
of the elements $C_w$, $w\klr\wnotl$, equals the right ideal
$C_{\wnotl}\hecke$.  Similarly
$\langle C_w\st w\kll\wnotl\rangle_A=
\hecke C_{\wnotl}$.
\end{proposition}
\begin{proof}
It is enough to prove the first equality,  the second being
a left analogue of the first.
The inclusion~$\supseteq$ follows immediately from the definition
of~$\klr$;  the inclusion~$\subseteq$ from~\cite[Lemma~2.11]{mp}.  
\end{proof}
\noindent
Setting $\xl:=\sum_{w\in\Wl} v_w T_w$ and
$\yl:=\sum_{w\in\Wl}\epsilon_w v_w^{-1}T_w$, we have, by~\cite[Theorem~1.1, Lemma~2.6\,(vi)]{kl}:
\begin{equation}\label{eq:xlyl}
   \xl=v_{\wnotl}C'_{\wnotl} \quad
\quad\quad\quad\quad\quad\quad
\yl=\epsilon_{\wnotl}v_{\wnotl}^{-1}C_{\wnotl}
\end{equation}
\begin{lemma}\mylabel{l:cfreemlin}
Let $\lambdad$ denote the partition~$(d+1,1,\ldots,1)$ of~$n$.
The two-sided ideal 
generated by $C_{\wnotlambdad}$ is a free $A$-submodule of~$\hecke$
with basis~$C_x$,
$\textup{\shape}(x)$ has more than~$d$ rows (or, equivalently,
$\textup{\shape}(x)\domleq\lambdad'$).
  \end{lemma}
  \begin{myproof}
Since $w_{0,\lambdad}$ has shape~$\lambdadp$ (see Proposition~\ref{p:pftinv}~(\ref{i:wnotl})),
it follows from the combinatorial
description of $\kllr$ in~\S\ref{sss:kllr} that $x\kllr w_{0,\lambdad}$ if 
and only if $\textup{\shape}(x)\domleq\lambdadp$.   So it is clear
from the definition of the relation~$\kllr$ (\S\ref{ss:klorders})
that the two-sided ideal $\hecke \cwnotld \hecke$ is contained in
$\langle C_x\st \textup{\shape}(x)\domleq\lambdadp\rangle_A$.
To show the reverse containment,  we first observe that $C_x$
belongs to the right ideal $\cwnotld\hecke$ 
in case $\mu:=\textup{\shape}(x)\domleq\lambdadp$
and $x=w_{0,\mu'}$, the longest element of its shape:
it is enough, by Proposition~\ref{p:pftinv2},  to show that $x\klr\wnotld$;
on the other hand,   by Proposition~\ref{p:pftinv}~(\ref{i:wnotldija}),  
$x\klr\wnotld$ is equivalent
to $x(1)>x(2)>\ldots>x(d+1)$,  which clearly holds for 
the elements~$x$ that we are considering.

Now suppose that~$x$ is a general element of \shape~$\mu\domleq\lambdadp$.
Proceed by induction on the domination order of~$\mu$.    
Let  $x\leftrightarrow
(P,Q)$ under RSK.   Then, on the one hand,
the association $C(P,Q)\leftrightarrow C(t_\mu,Q)$
gives an $\hecke$-isomorphism between the right cell modules~$\rcmod(x)$ and
$\rcmod{(v)}$, where $v$ is the permutation corresponding under RSK
to~$(t_\mu,Q)$ (\S\ref{ss:cellmod}); on the other,     
since~$v\leftrightarrow(t_\mu,Q)$ is right
equivalent to~$w_{0,\mu'}\leftrightarrow(t_\mu,t_\mu)$,  there exists,
by~Proposition~\ref{p:pftinv2},   an element~$h$ in~$\hecke$ such that
$C(v)=C_{w_{0,\mu'}}h$; so that, by the definition of right cell modules,
\[
C_x\equiv C_u h\pmod{\langle C_y \st y\klrlneq u\rangle_A}
\]
where $u\leftrightarrow(P,t_\mu)$ under RSK.
Now, $y\klrlneq u$ implies,
by~\S\ref{sss:unrelate},~\ref{sss:kllr},  
$\textup{\shape}(y)\triangleleft\mu\domleq\lambdadp$;  and, by the induction
hypothesis, $C_y\in\hecke \cwnotld\hecke$.  
As to $C(P,t_\mu)$, being left equivalent to
$C(t_\mu,t_\mu)$, it belongs, once again by Proposition~\ref{p:pftinv2},
to the left ideal~$\hecke C_{w_{0,\mu'}}$,  which
as shown in the previous paragraph is contained in~$\hecke \cwnotld \hecke$.
Thus $C_x=C(P,Q)\in\hecke\cwnotld\hecke$,  and we are done.
  \end{myproof}

\mysubsubsection{Proof of Theorem~\ref{t:invariant} given
Lemma~\ref{l:cfreemlin}}\mylabel{sss:pftinv}
\newcommand\jlift{\tilde{J}}
As seen in~\S\ref{sss:hspecial},  $\base\symm_n$ is the specialization
of the Hecke algebra~$\hecke$:  $\base\symm_n\simeq\hecke_\base:=
\hecke\tensor_A\base$, where $\base$ is an $A$-algebra via the
natural ring homomorphism $A\to\base$ defined by $v\mapsto1$.
Under the map $\hecke\to\hecke\tensor_A\base$ given by $x\mapsto x\tensor1$,
the image of $C_{\wnotld}$ is $C_{\wnotld}\tensor 1=y_d$, by 
Eq.~(\ref{eq:xlyl}).   
Denoting by~$\jlift$
the two-sided ideal of~$\hecke$ generated by~$C_{\wnotld}$,  we thus have
$\hecke/\jlift\tensor_A\field\simeq \base\symm_n/\jnd$. 

On the other hand, combining Lemma~\ref{l:cfreemlin} with the observation
in~\S\ref{sss:observe},   we see that $\hecke/\jlift$ is a free $A$-module
with basis~$T_x$,  as $x$ varies over permutations of whose \shape{s} have
at most~$d$ rows.    The image of $T_x$ in~$\base\symm_n/\jnd$ being the
residue class of the corresponding permutation~$x$,   the theorem
is proved.\hfill$\Box$
\mysubsection{A `monomial' basis for the $\glv$-invariant sub-algebra of the tensor algebra of~$\EndV$}%
\mylabel{ss:tensor}  Let $\field$ be a commutative ring with identity such
that no non-zero polynomial in one variable
over~$\field$ vanishes identically as a function on~$\field$.
Let $V$ be a free module over~$\field$ of finite rank~$d$.
Let $T:=T(\EndV)$ denote the tensor algebra $\oplus_{n\geq 0}
(\EndV)^{\tensor n}$.  The action of the group~$\glv$ of units in~$\EndV$
on $T$ preserves the algebra
structure, so the ring~$T^{\glv}$ of $\glv$-invariants is a sub-algebra.   
It also preserves degrees, so \[T^{\glv}\simeq
\oplus_{n\geq0}((\EndV)^{\tensor n})^{\glv}=\oplus_{n\geq0}\Endo_{\glv}(V^{\tensor n}).\]
By the classical theorem quoted in~\S\ref{ss:iinv} from~\cite{dp},
we have, for every~$n\geq0$, an isomorphism of $\field$-algebras 
\[\Theta_n\!:\ \ \field\symm_n/\jnd\simeq\Endo_{\glv}(V^{\tensor n})\]
where $\field\symm_n$ is the group algebra of the symmetric group~$\symm_n$,
and $\jnd$ the two sided ideal as defined 
in the statement of the quoted theorem.

\newcommand\basis{\mathfrak{B}}
Now, Theorem~\ref{t:invariant} gives us a $\field$-basis for 
$\field\symm_n/\jnd$.   Taking the image under~$\Theta_n$ gives a
basis for~$\Endo_{\glv}(V^{\tensor n})$.
Taking the disjoint union over~$n$ of these bases
gives a basis---call it~$\basis$---for~$T^{\glv}$, which 
has an interesting property---see the theorem below---which 
explains the appearance of term `monomial'
in the title of this subsection.
\begin{theorem}
  \mylabel{t:invapp1}
The basis~$\basis$ of~$T^{\glv}$ defined above is closed under products.
\end{theorem}
\begin{proof}
\newcommand\subalg{\mathfrak{P}}
\newcommand\subalgn{\subalg_n}
In fact, we get a description of the $\field$-algebra $T^{\glv}$ as follows.
Consider the space 
$\mathfrak{S}:=\oplus_{n\geq 0}\field\symm_n$ with the following
multiplication:  for $\pi$ in $\symm_m$ and $\sigma$ in~$\symm_n$,
$\pi\cdot\sigma$ is the permutation in $\symm_{m+n}$ that, as a self-map
of $[m+n]$, is given by
\[ \pi\cdot\sigma(i):=\left\{
  \begin{array}{ll}
    \pi(i) & \textup{if $i\leq m$}\\
    \sigma(i-m)+m & \textup{if $i\geq m+1$}
  \end{array}
  \right.\]
For each $n$, consider the subspace~$\subalgn$ of $\field\symm_n$ spanned by 
permutations that have no decreasing sub-sequence of length more than~$d$.
The direct sum $\subalg:=\oplus_{n\geq0}\subalgn$ is a sub-algebra 
of~$\mathfrak{S}$.

The restriction to~$\subalgn$ of the canonical map $\field\symm_n
\to \field\symm_n/\jnd$ is a vector space isomorphism 
(Theorem~\ref{t:invariant}).   Thus 
$\oplus_{n\geq0}\Theta_n$ is a vector space isomorphism of the 
algebra~$\subalg$ onto~$T^{\glv}$.   It is evidently also
an algebra isomorphism.
\end{proof}
\mysubsection{Application to rings of polynomial invariants}%
\mylabel{ss:poly}
In this subsection, $\field$ denotes a field of characteristic~$0$
and $V$ a $\field$-vector space of finite dimension~$d$.
Consider the ring of $\textrm{GL}_\field(V)$-invariant polynomial functions 
on $(\End_\field{V})^{\times m}$.    It is a direct sum of homogeneous invariant
functions,   for the action preserves degrees.
It is spanned, as a vector space, by
products of traces in words (see for example \cite[\S1]{procesi},~\cite{razmyslov}).
More precisely, for a fixed degree~$n$, given a permutation
$\sigma$ of~$n$ elements and a map $\nu$ of $[n]$ to $[m]$,
consider the function~$f(\sigma,\nu)$ defined as follows:
writing $\sigma$ as a product $(i_1i_2\cdots)
(i_{k+1}i_{k+2}\cdots)\cdots(i_{p+1}i_{p+2}\cdots)$ of disjoint cycles,
\begin{equation*}\begin{split}
f(\sigma,\nu):=\textup{Trace}(A_{\nu({i_1})}A_{\nu({i_2})}\cdots)
\textup{Trace}(A_{\nu{(i_{k+1})}}A_{\nu{(i_{k+2})}}\cdots)\cdots\\
\quad\cdots
\textup{Trace}(A_{\nu{(i_{p+1})}}A_{\nu({i_{p+2})}}\cdots)
\end{split}\end{equation*}
As $\sigma$ and $\nu$ vary, the $f(\sigma,\nu)$ span the
space of invariants of degree~$n$.

This fact is proved by observing that every polynomial
invariant arises as the specialization of a multilinear invariant
(by the restitution process). Thus,
thanks to Theorem~\ref{t:invariant}, we can restrict the permutation $\sigma$
to have no decreasing sub-sequence of length more than~$d$,  and still
the $f(\sigma,\nu)$ would span.   We state this formally:
\begin{theorem}
  \mylabel{t:poly1}
The invariant functions~$f(\sigma,\nu)$, as $\sigma$ varies over 
permutations that do not have any decreasing sub-sequence of 
length exceeding~$d$, form a $\field$-linear spanning set for the
ring of $\textrm{GL}_\field{V}$-invariant polynomial functions on
$(\End_\field{V})^{\times n}$.
\end{theorem}
\mysubsubsection{Picture invariants}\mylabel{sss:picture}
Set $V_b^t:={V^*}^{\tensor b}\tensor V^{\tensor t}$ and consider the
ring of polynomial $\textrm{GL}_\field(V)$-invariant functions on the space
$V_{b_1}^{t_1}\times\cdots\times V_{b_s}^{t_s}$ of several tensors.     
In~\cite[\S3]{dks},
the notion of a `picture invariant' is introduced, generalizing
the functions~$f(\sigma,\nu)$ defined above.   Picture invariants
span the space of invariant polynomial functions~(\cite[Proposition~7]{dks}).
Just as in the special case of $(\EndV)^{\times n}$ discussed above,
thanks to Theorem~\ref{t:invariant}, we have:
\begin{theorem}
  \mylabel{t:poly2}
Only those picture invariants with underlying permutations 
having no decreasing sub-sequences of length exceeding~$d$ suffice to span 
as a $\field$-vector space the ring of $\textrm{GL}_\field(V)$-invariant
polynomial functions on the space $V_{b_1}^{t_1}\times\cdots\times V_{b_s}^{t_s}$
of several tensors.
\end{theorem}

%% file: tabloidproof.tex
\mysection{Proof of Theorem~\ref{t:tabloid}}
\label{s:part3}\mylabel{s:tabloidproof}
\noindent
In this section,  we first
prove Theorem~\ref{t:tabloid}~(\S\ref{ss:itabloid}).   
We then show that it holds also over the integers and fields
of characteristic~$0$ (\S\ref{ss:tabloidz}),  but not in general
over a field of positive characteristic (Example~\ref{eg:tabloid}).
For comments on its Hecke analogue, see~\S\ref{ss:ipf}.

As pointed out in~\S\ref{ss:itabloid},  the results of the recent
paper~\cite{dn} are related to Theorem~\ref{t:tabloid} and its
Hecke analogue.
\mysubsection{Proof of Theorem~\ref{t:tabloid}}\mylabel{ss:tabloidpf}
Let $n$ be a positive
integer,  $\lambda$ a partition of~$n$, and $\repl:\fieldc\symm_n\to
\End_\fieldc{\fieldc\tabloidl}$ the map defining the representation of
$\symm_n$ on tabloids of shape~$\lambda$ (\S\ref{ss:tabloid}).
The proof follows by combining 
the observations in~\S\ref{ss:projections} and~\S\ref{sss:observe}
with the following two facts:  
\begin{enumerate}
\item the decomposition into irreducibles of the 
representation~$\fieldc\tabloidl$
is given by $\fieldc\tabloidl=\oplus_{\mu\domgeq\lambda}
(\specht^\mu_\fieldc)^{m(\mu)}$,  where~$\specht^\mu_\fieldc$ are the Specht
modules (\S\ref{ss:specht}),  $\domgeq$ is the domination relation on
partitions (\S\ref{ss:domination}), and the multiplicities $m(\mu)$ are
positive.\footnote{In fact, $m(\mu)$ is the number of
`semi-standard tableaux of shape~$\mu$ and content~$\lambda$'.}
\item the Specht module $\specht^{\mu}_\fieldc$ 
is $\symm_n$-isomorphic to the right cell module~$\rcmod(\mu)_\fieldc$
(defined in~\S\ref{ss:cellmod});
\end{enumerate}
Both facts are well known.
For~(1), see for example~\cite[Theorem~2.11.2, Corollary~2.4.7]{sagan}. 
For~(2), we could refer to~\cite{gm} or~\cite{naruse}.   But in fact
we will recall in some detail in~\S\ref{ss:mpiso}  
the following more general fact from~\cite{dj,mp}:  %
Specht modules can be defined over the Hecke algebra~$\hecke$
and are isomorphic to the corresponding right cell modules.

Since the multiplicities~$m(\mu)$ in~(1) above are positive,
the kernel of~$\repl$ is the same as that of the map~$\repl':\fieldc\symm_n
\to\End_\fieldc(\oplus_{\mu\domgeq\lambda}\specht^{\mu}_\fieldc)$.   The image
of~$\repl'$ is clearly contained 
in~$\oplus_{\mu\domgeq\lambda}\End_\fieldc\specht^{\mu}_\fieldc$. Since
the $\specht^{\mu}_\fieldc$ are non-isomorphic for distinct~$\mu$,%
\footnote{\label{f:three}
This is well-known.  It also follows from the isomorphism in~(2) and
the corresponding fact for cell modules proved in~\S\ref{ss:wstcon}.}
it follows from a density argument (see for example~\cite[Chapter~8, \S4,
No.~3, Corollaire~2]{bouralg}) that $\repl'$ maps 
onto~$\oplus_{\mu\domgeq\lambda}\End_\fieldc\specht^{\mu}_\fieldc$. 
Since $\dim\specht^{\mu}_\fieldc=d(\mu)$, where $d(\mu)$ is the number
of standard tableaux of shape~$\mu$, and
the $\specht^{\mu}_\fieldc$ as $\mu$ varies over all partitions of~$n$
are a complete set of irreducible representations,%
\footnote{Same comment as in footnote~\ref{f:three} applies to both assertions.}
we obtain, by counting dimensions:
\[\dim{\ker{\repl'}}=\dim{\fieldc\symm_n}-\dim{(
\oplus_{\mu\domgeq\lambda}\End_\fieldc\specht^{\mu}_\fieldc)}
= \sum_{\mu\partition n}d(\mu)^2
- \sum_{\mu\domgeq\lambda} d(\mu)^2\\
=
\sum_{\mu\not\domgeq\lambda} d(\mu)^2
\]

Now consider $\fieldc\symm_n$ as the specialization of the Hecke algebra
$\hecke$ as follows (\S\ref{sss:hspecial}):  $\fieldc\symm_n\simeq
\hecke\tensor_A\fieldc$,  where $\fieldc$ is an $A$-algebra via the
map~$A\to\fieldc$ defined by~$v\mapsto1$.
By the observation~\S\ref{ss:projections}, the
images $C_w\tensor1$ in~$\hecke\tensor_A\fieldc\simeq\fieldc\symm_n$ 
of the Kazhdan-Lusztig basis elements~$C_w$ of~$\hecke$ (\S\ref{ss:klbasis})
belong to the kernel of~${\repl'}$
if~$\textup{RSK-shape}(w)\not\domgeq\lambda$.   The number of such
$w$ being equal to~$\sum_{\mu\not\domgeq\lambda}d(\mu)^2$,  which as observed above 
equals $\dim{\ker{\repl'}}$,   we conclude that 
\begin{equation}
  \label{eq:kerrepl}
\ker{\repl}=   
\ker{\repl'}=
\langle C_w\tensor1\st \textup{\shape}(w)\not\domgeq\lambda\rangle_\fieldc.
\end{equation}

By observation~\S\ref{sss:observe}, the images of $T_w\tensor1$,
$\textrm{\shape}(w)\domgeq\lambda$,  form a basis for 
$\hecke\tensor_A\fieldc/\langle C_x\tensor 1\st \textup{\shape}(w)\not\domgeq
\lambda\rangle_\fieldc\simeq \fieldc\symm_n/\ker{\repl'}$.   But the image
in~$\fieldc\symm_n$ of $T_w\tensor1$ is the permutation~$w$.
This completes the proof of Theorem~\ref{t:tabloid}.\hfill$\Box$
\mysubsection{$\jnd$ in~$\fieldc\symm_n$ equals
$\ker{\replnd}$}  \mylabel{sss:jndkerreplnd}
We now justify the claim made in~\S\ref{ss:itabloid} that the
ideal $\jnd$ in~$\fieldc\symm_n$ equals $\ker{\replnd}$.   On the one hand,
as is easily seen, the generator~$y_d$ of the two sided ideal~$\jnd$
belongs to~$\ker{\replnd}$.   Indeed, given a tabloid~$\{T\}$ of shape~$\lnd$,
there evidently exist integers $a$ and $b$, with $1\leq a,b\leq d+1$, 
that appear in the same row of~$T$.   This implies that the transposition
$(a,b)$ fixes~$\{T\}$.  Writing $\symm_{d+1}$ as a disjoint union 
$S\cup S(a,b)$ (for a suitable choice of a subset~$S$), we have
$y_d\{T\}=\sum_{\sigma\in\symm_{d+1}}\sign{(\sigma)}\sigma\{T\}
=\sum_{\sigma\in S}\sign(\sigma)(\sigma-\sigma(a,b))\{T\}=0$.

On the other hand, as computed in the proof of Theorem~\ref{t:tabloid}
above, $\ker{\rep_{\lnd}}$ as a $\fieldc$-vector space has dimension
$\sum_{\mu\not\domgeq\lnd}d(\mu)^2$.   
It suffices therefore to show that $\jnd$ too has this same dimension.
It follows from Lemma~\ref{l:cfreemlin} 
that $\jnd$ has dimension $\sum_{\mu\domleq\lambdadp}d(\mu)^2$, where $\lambdad$
is the partition of $(d+1,1,\ldots,1)$ of~$n$ and $\lambdadp$ denotes
its transpose (\S\ref{sss:invprelims}).   But $\mu\not\domgeq\lnd$
if and only if $\mu$ has more than~$d$ rows if and only if
$\mu\domleq\lambdadp$.\hfill$\Box$
\mysubsection{Theorem~\ref{t:tabloid} holds over the integers and fields of
characteristic~$0$}%
\mylabel{ss:tabloidz}
\newcommand\replz{\rep_{\lambda,\integers}}
We first argue that Theorem~2 holds with $\integers$ coefficients in place
of $\fieldc$ coefficients.   Let $\replz$ be the map
$\integers\symm_n\to\End_\integers{\integers\tabloidl}$ defining the tabloid
representation.   We claim that Eq.~(\ref{eq:kerrepl}) holds over~$\integers$:
\begin{equation}
  \label{eq:kerreplz}
  \ker{\replz}=\langle C_w\tensor1\st \textup{\shape}(w)\not\domgeq\lambda
\rangle_\integers
\end{equation}
Once this is proved,  the rest of the argument is the same as in the complex
case: namely, use observation~\S\ref{sss:observe}.

We first show the containment~$\supseteq$.
We have $(C_w\tensor1)\fieldc\tabloidl=(C_w\tensor1)\integers\tabloidl\tensor_\integers\fieldc$ (by flatness of~$\fieldc$ over~$\integers$).    
Since $(C_w\tensor1)\integers\tabloidl$ is a submodule of the free module
$\integers\tabloidl$, it is free.   By Eq.~(\ref{eq:kerrepl}),
$(C_w\tensor1)\fieldc\tabloidl=0$ if $\textup{\shape}(w)\not\domgeq\lambda$,
so~$\supseteq$ holds.

\newcommand\module{\frak{m}}
To show the other containment, 
set $\module=
\langle C_w\tensor1\st \textup{\shape}(w)\not\domgeq\lambda\rangle_\integers$,
and
consider~$\ker{\replz}/\module$.
Since $\integers\symm_n/\module$ is free, so is its submodule
$\ker{\replz}/\module$,
and we have \[\frac{\ker{\replz}}{\module}\tensor_\integers\fieldc
=\frac{\ker{\replz}\tensor_\integers\fieldc}{\module\tensor_\integers\fieldc}
=\frac{\ker{\replz}\tensor_\integers\fieldc}{\langle C_w\tensor1\st \textup{\shape}(w)\not\domgeq\lambda\rangle_\fieldc}.\]
By the flatness of~$\fieldc$ over~$\integers$,  
we have $\ker{\replz}\tensor_\integers\fieldc=\ker{\repl}$.   
The last term in the above display vanishes
by Eq.~(\ref{eq:kerrepl}), and so~$\subseteq$ holds (since 
$\ker{\replz}/\module$ is free).   The proof of Theorem~\ref{t:tabloid}
over~$\integers$ is complete.

Let now $\field$ be a field of characteristic~$0$ and $\rep_{\lambda,\field}$
the map $\field\symm_n\to\End_\field\field\tabloidl$ defining the representation
on tabloids of shape~$\lambda$.   The analogue of
Eqs.~(\ref{eq:kerrepl}) and (\ref{eq:kerreplz}) holds over~$\field$,
since, by the flatness of~$\field$ over~$\integers$,  we have
$\ker{\rep_{\lambda,\field}}=\ker{\replz}\tensor_\integers\field$. 
Now use observation~\S\ref{sss:observe} as in the earlier cases to finish
the proof of Theorem~\ref{t:tabloid} over~$\field$.\hfill$\Box$
\begin{example}
  \mylabel{eg:tabloid}  \begin{upshape}
Theorem~\ref{t:tabloid} does not hold in general over a field~$\field$ of
positive characteristic.
We give an example of a non-trivial linear
combination of permutations of 
\shape\ dominating~$\lambda$ that acts trivially on the tabloid representation 
space~$\field\tabloidl$.
Let~$\field$ be a field of characteristic~$2$. Let $n=4$ and $\lambda=(2,2)$.
Let us denote a permutation in~$\symm_4$ by writing down in
sequence the images under it of $1$ through~$4$: e.g., $1243$ denotes the
permutation~$\sigma$ defined by
$1\sigma=1$, $2\sigma=2$,  $3\sigma=4$, and $4\sigma=3$.  It is
readily seen that the eight permutations in the display 
below are all of shape $(3,1)$ and that their sum 
acts trivially on~$\field\tabloidl$.   
\[2134, \quad
2341, \quad
2314, \quad
1342,\quad
3124,\quad
1243, \quad
4123, \quad
1423.\]
\end{upshape}\end{example}

%% file: proof.tex
\mysection{Proof of Theorem~\ref{t:irrep}}%
\mylabel{s:proof}
\noindent
In this section, $\base$ denotes a field, $a$ an invertible
element of~$\base$,  and notation is fixed as in~\S\ref{s:setup}.
It is assumed throughout this section (but not in the later ones)
that~$\hecke_\base$ is semisimple (see~\S\ref{sss:hspecial}).
Combining the structure theory of semisimple algebras
with the results 
in~\S\ref{ss:cellsandkrs}--\ref{ss:projections},  we derive 
Theorem~\ref{t:irrepgen} as a consequence.
Theorem~\ref{t:irrep} of~\S\ref{s:intro} is the special case
of this theorem when $\base$ is the complex field and $a=1$. 
As mentioned in~\S\ref{ss:ipf},  the proof in effect says that 
Theorem~\ref{t:klbasis} is a consequence of the cellularity in
the sense of~\cite{glinv} of the Kazhdan-Lusztig basis.

But first,  we need to recall the structure theory,
which we do in~\S\ref{ss:wburn} in a form suited to our
context,  and then the proof of irreducibility of the
(right) cell modules,
which we do in~\S\ref{ss:wstcon} following the argument in~\cite[\S5]{kl}.   
\mysubsection{A recap of the structure theory of semisimple algebras}%
\mylabel{ss:wburn}
\newcommand\algebraa{\mathfrak{A}}
\newcommand\divv{E_V}
\newcommand\dimvdivv{n_V}
\newcommand\divvopp{D_V}
\newcommand\matrixm{\mathcal{M}}
The facts recalled here are all well known: see, e.g.,~\cite{bouralg},
\cite[page~218]{gp}.
Let~$V$ be a simple (right) module for a semisimple algebra~$\algebraa$
of finite dimension over~$\base$.
Then the endomorphism ring~$\End_{\algebraa}V$ is a division algebra 
(Schur's Lemma),  say $\divv$.  Being a subalgebra of~$\End_\field V$,
it is finite dimensional as a vector space over~$\field$,
and~$V$ is a finite dimensional vector space over it.  
Set $\dimvdivv:=\dim_{\divv}V$.
The ring~$\End_{\divv}V$ of endomorphisms of~$V$ as a $\divv$-vector space
can be identified (non-canonically, depending upon a choice of basis)
with the ring $\matrixm_{\dimvdivv}(\divvopp)$ of matrices of size 
$\dimvdivv\times\dimvdivv$ with entries in the opposite algebra~$\divvopp$
of~$\divv$.    The natural ring homomorphism~$\algebraa\to\Endo_{\divv}V$
is a surjection (density theorem).    

There is an isomorphism of algebras (Wedderburn's structure theorem):
\begin{equation}
  \label{eq:wburn}
\algebraa\simeq \prod_V \Endo_{\divv}V\simeq\prod_V \matrixm_{\dimvdivv}(\divvopp), 
\end{equation}
where the product is taken over all (isomorphism classes of) simple modules.
There is a single isomorphism class of simple modules for the simple 
algebra~$\Endo_{\divv}V$,  namely that of $V$ itself, 
and its multiplicity is~$\dimvdivv$ in a direct sum decomposition into simples
of the right regular representation of~$\Endo_{\divv}V$.   Thus~$\dimvdivv$ is
also the multiplicity of~$V$ in the right
regular representation of~$\algebraa$.    And of course
\begin{equation}
  \label{eq:dimandmult}
\dim_\field V
=\dimvdivv(\dim_\field\divv)\geq\dimvdivv  
\end{equation}

The hypothesis of the following proposition admittedly appears
contrived at first sight,
but it will soon be apparent (in~\S\ref{ss:klbasis}) that it is tailor-made
for our situation.
\begin{proposition}
  \mylabel{p:wburn}
Let $W_1$, \ldots, $W_s$ be $\algebraa$-modules of respective dimensions
$d_1$, \ldots, $d_s$ over~$\field$.     Suppose that the right regular
representation of~$\algebraa$ has a filtration in which the quotients
are precisely $W_1^{\oplus d_1}$, \ldots, $W_s^{\oplus d_s}$.    Then
\begin{enumerate}
\item $\Endo_{\algebraa}W_i=\field$ and $W_i$ is absolutely irreducible,
$\forall$ $i$, $1\leq i\leq s$.
\item $W_i$ is not isomorphic to~$W_j$ for $i\neq j$.
\item $W_i$, $1\leq i\leq s$, are a complete set of simple $\algebraa$-modules.
\item $\algebraa\simeq\prod_{i=1}^{s}\Endo_\field{W_i}$.
\end{enumerate}
\end{proposition}
\begin{proof}
  Let $V$ be a simple submodule of~$W_i$.   Then $\dim_\field V\leq d_i$.
The hypothesis about the filtration implies that the multiplicity of~$V$
in the right regular representation is at least~$d_i$.  
From Eq.~(\ref{eq:dimandmult}), we conclude that $d_i=\dimvdivv$ and
$\dim_\field\divv=1$.    So $V=W_i$ is simple and $\divv=\field$.
\newcommand\kbar{{\bar{\field}}}
If $\kbar$ denotes an algebraic closure of~$\field$,  then 
\[\Endo_{\algebraa\tensor_{\field}\kbar}(V\tensor_\field\kbar)=(\Endo_{\algebraa}{V})\tensor_\field\kbar=\field\tensor_\field\kbar=\kbar.\]
So $V$ is absolutely irreducible and~(1) is proved.

If $W_i\simeq W_j$ for $i\neq j$,  then the multiplicity of~$W_i$ in the
right regular representation would exceed~$d_i$ contradicting 
Eq.~(\ref{eq:dimandmult}).    This proves~(2).   Since every simple
module has positive multiplicity in the right regular representation,
(3) is clear.  Finally, (4) follows from (1) and Eq.~(\ref{eq:wburn}).
\end{proof}
\mysubsection{Irreducibility and other properties of the cell modules}%
\label{s:irreducible}\mylabel{ss:wstcon}   The proof of 
Theorem~\ref{t:klwstcon} below follows~\cite[\S5]{kl}.    We give it in detail
here for the proof in~\cite{kl} seems sketchy.    
\begin{theorem}
  \mylabel{t:klwstcon} \textup{(\cite[\S5]{kl})}
Assume that $\heckek:=\hecke\tensor_A\base$ is semisimple.  Then
\begin{enumerate}
\item $\Endo_{\heckek}\rcmodlk=\field$ and $\rcmodlk$ is absolutely irreducible,
for all $\lambda\partition n$.
\item $\rcmodlk\not\simeq \rcmodmk$ for partitions $\lambda\neq\mu$
of $n$.
\item $\rcmodlk$, $\lambda\partition n$, are a complete set of 
simple $\heckek$-modules.
\item $\heckek\simeq\prod_{\lambda\partition n}\Endo_\field{\rcmodlk}$.
\end{enumerate}
\end{theorem}
\begin{proof}
\newcommand\filter{\mathfrak{F}}
  By Proposition~\ref{p:wburn},  it is enough to exhibit a filtration
of the right regular representation of $\heckek$ in which the quotients
are precisely $\rcmodlk^{\oplus\dl}$, $\lambda\partition n$, each occurring
once.   We will in fact exhibit a
decreasing filtration $\filter=\{F_i\}$ by right ideals (in fact,
two sided ideals) of~$\hecke$ in which the quotients $F_i/F_{i+1}$ are
precisely $\rcmodl^{\oplus\dl}$, $\lambda\partition n$, each occurring
once.   Since $\rcmodl$ are free $A$-modules,  it will follow that 
$\filter\tensor_A\base$ is a filtration of~$\heckek$ whose quotients
are $\rcmodlk^{\oplus\dl}$,  and the proof will be done.

Let $\succeq$ be a total order on partitions of $n$ that refines the
dominance partial order~$\domgeq$.   Let $\lambda_1\succ\lambda_2\succ\ldots$
be the full list of partitions arranged in decreasing order with respect
to $\succeq$.    Set $F_i:=\langle C_w\st \textrm{\shape}(w)\preceq\lambda_i\rangle_A$.
It is enough to prove the following:
\begin{enumerate}
\item The $F_i$ are right ideals in~$\hecke$ (they are in fact
two sided ideals).
\item $F_i/F_{i+1}\simeq \rcmod(\lambda_i)^{\oplus d(\lambda_i)}$.
\end{enumerate}

It follows from the definition in~\S\ref{ss:klorders} of the relation
$\kllr$ that,
for any fixed permutation~$w$, $\langle C_x\st x\kllr w\rangle_A$ is a
two sided ideal of~$\hecke$.  
But $x\kllr w$ if and only if $\rskshape{(x)}\domleq\rskshape{(w)}$,
by the characterization in~\S\ref{sss:kllr}.
Thus,  $\langle C_x\st \rskshape{(x)}\domleq \lambda\rangle_A$
is a two sided ideal, and $F_i$ being equal to the sum
$\sum_{j\geq i}\langle C_x\st \rskshape{(x)}\domleq\lambda_j\rangle_A$
of two sided ideals is a two sided ideal.   This proves~(1).

To prove~(2),    let $S_1$, $S_2$, \ldots\ be the distinct
right cells contained in the two sided cell corresponding to 
shape~$\lambda_i$. It follows from the assertions in~\S\ref{sss:cells} that
there are $d(\lambda_i)$ of them and the cardinality of each is
$d(\lambda_i)$.    Fix a permutation $w$ of shape $\lambda_i$.
Consider the right cell module $\rcmod(w)$, which by definition
is the quotient of the right ideal $\langle C_x\st x\klr w\rangle_A$
by the right ideal $\langle C_x\st x\klrlneq w\rangle_A$.  If $x\klr w$
then evidently $x\kllr w$ and (by \S\ref{sss:kllr}) $\rskshape{(x)}\domleq
\lambda_i$,  so $\rskshape{(x)}\preceq\lambda_i$.   Thus we have
a map induced by the inclusion: $\langle C_x\st x\klr w\rangle_A
\to F_i/F_{i+1}$.

We claim that the above map descends to an injective map from the
quotient~$\rcmod(\lambda_i)$.   It descends because $x\klrlneq w$
implies $x\kllrlneq w$:   if $x\kllrsim w$,  then $x\klrsim w$
by~\S\ref{sss:unrelate}.  To prove that the map from~$\rcmod(\lambda_i)$
is an injection,  let $\sum_{x\kllr w} a_xC_x$ belong to $F_{i+1}$ with
$a_x\in A$.   Suppose that $a_x\neq 0$ for some fixed~$x$.
Then,  since the $C_y$ form an $A$-basis of~$\hecke$,
we conclude that $\rskshape{(x)}\preceq\lambda_{i+1}$,
so $\rskshape{(x)}\neq\lambda_i$,  and (by~\S\ref{sss:kllr})
$x\kllrlneq w$.   But this means $x{\not\sim}_R w$,
so $x\klrlneq w$, and thus the image in $\rcmod{(\lambda_i)}$
of $\sum_{x\kllr w}a_xC_x$ vanishes.   

The image of $\rcmod(w)$ in $F_i/F_{i+1}$ is spanned by the classes
$\bar{C_x}$, $x\klrsim w$.   Choosing $w_1$ in $S_1$, $w_2$ in $S_2$, \ldots\ 
we see that the images of $\rcmod(w_1)$, $\rcmod(w_2)$, \ldots\ in
$F_i/F_{i+1}$ form a direct sum (for the~$C_x$ are an $A$-basis of~$\hecke$).
The $\rcmod(w_j)$ are all isomorphic to~$\rcmodl$ (see~\S\ref{ss:cellmod}).
This completes the proof of (2) and also of the theorem.
\end{proof}
\mysubsection{Kazhdan-Lusztig basis in endomorphisms
of modules}\mylabel{ss:klendo}
\begin{theorem}\mylabel{p:klbasis}\mylabel{t:klbasis}\label{t:irrepgen}
Assume that $\heckek:=\hecke\tensor_A\field$ is semisimple.
For $\lambda$ a partition of~$n$, the images in $\Endo{\rcmodlk}$ of
the Kazhdan-Lusztig basis elements
$C_x$, $\textup{\shape}(x)=\lambda$, 
form a basis (for $\Endo{\rcmodlk}$).\end{theorem}
\begin{myproof}
By Theorem~\ref{t:klwstcon}~(4),  $\heckek\simeq\oplus_{\lambda\vdash n}
\Endo{\rcmodlk}$.
The projections to $\Endo{\rcmod(\nu)_\field}$ of $C_x$,
$\textup{\shape}(x)\domleq\lambda$,  vanish if $\nu\not\domleq\lambda$
(\S\ref{ss:projections}).   Therefore the projections of the same elements
to $\oplus_{\mu\domleq\lambda}\Endo{\rcmod(\mu)_\field}$ form a basis: note that the 
number of such elements equals 
$\sum_{\mu\domleq\lambda}\dim{\Endo_{\rcmod(\mu)_\field}}$.
Again by~\S\ref{ss:projections}, the projections
of $C_x$, $\textup{\shape}(x)\doml\lambda$, 
vanish in $\Endo{\rcmod(\lambda)_\field}$.   This implies that
the projections of $C_x$, $\rskshape{(x)}=\lambda$, in~$\Endo{\rcmodlk}$
form a spanning set.   Since the number of such $C_x$ equals 
$\dim{\Endo{\rcmodlk}}$,
the theorem follows.\end{myproof}
\begin{theorem}\mylabel{p:p2kl}\mylabel{t:p2kl}
Assume that $\heckek:=\hecke\tensor_A\field$ is semisimple.
Let $U$ be a finite dimensional representation
of~$\hecke_\field$ and $\subsets$ the subset of partitions~$\lambda$ of~$n$
such that $\rcmodlk$ appears in a decomposition of~$U$ into irreducibles.
Then the images in $\Endo{U}$ of $C_x$, $x\in\symm_n$ such that
$\textup{\shape}(x)\in\subsets$, form a basis for the image of~$\heckek$
(under the map $\heckek\to\Endo{U}$ defining~$U$).
\end{theorem}
\begin{myproof}
It is enough to prove the assertion assuming 
$U=\oplus_{\lambda\in\subsets}\rcmodlk$.   The image of $\heckek$ in
$\Endo{U}$ is $\oplus_{\lambda\in\subsets}\Endo{\rcmodlk}$ (Theorem~\ref{t:klwstcon}~(1),  density theorem, and~\cite[Corollaire~2, page~39]{bouralg}).
Proceed by induction
on the cardinality of~$\subsets$.    It is enough to show that the
relevant images in $\Endo{U}$ are linearly independent,   for 
their number equals the dimension 
of~$\oplus_{\lambda\in\subsets}\Endo{\rcmodlk}$.
Suppose that a linear combination of the images vanishes.   
Choose~$\lambda\in\subsets$ such that
there is no~$\mu$ in $S$ with $\lambda\triangleleft\mu$.    
Projections to~$\Endo{\rcmodlk}$ of all $C_x$, $\lambda\neq
\textup{\shape}(x)\in\subsets$, vanish (\S\ref{ss:projections}).
So projecting the linear combination to $\Endo{\rcmodlk}$ and using
Theorem~\ref{p:klbasis},  we conclude that the coefficients of
$C_x$, $\textup{\shape}(x)=\lambda$, are all zero.  
The induction hypothesis applied to~$\subsets\setminus\{\lambda\}$
now finishes the proof.\end{myproof}
\bexample\mylabel{x:counter}
The purpose of this example is to show that images
in~$\Endo{\rcmodl_\fieldc}$ of permutations
of \shape\ $\lambda$ do not in general form a basis of $\Endo{\rcmodl_\fieldc}$.
Let $n=4$ and $\lambda=(2,2)$.     Then $\rcmodl_\fieldc$ is the unique
$2$~dimensional complex irreducible representation of~$\symm_4$.  Consider
the action of~$\symm_4$ on partitions of $\{1,2,3,4\}$ into two sets
of two elements each.   There being three such partitions,  we
get a map $\symm_4\to\symm_3$,  which is surjective and has kernel
$\{\textup{identity},
(12)(34),\,(13)(24),\,(14)(23)\}$.   
Pulling back the $2$-dimensional complex irreducible representation of~$\symm_3$
via the above map,  we get~$\rcmodl_\fieldc$.
The permutations of shape $\lambda$ are 
$(13)(24)$, $(1342)$, $(1243)$, and $(12)(34)$.
The first and last of these act as identity on~$\rcmodl_\fieldc$.
\eexample

%% file: matrixg.tex
\mysection{The matrix $\matrixgl$ and a formula for its determinant}%
\mylabel{s:matrixgl} 
\noindent
Let $\lambda$ be a fixed partition of $n$. Our goal in this section
is to study the action of
the elements~$C_w$, $\textup{\shape}(w)=\lambda$, on the right cell
module~$\rcmodl$.   More specifically, it is to state 
Theorem~\ref{t:formula}.    The motivation for this was already 
indicated (see~\S\ref{ss:iposchar}): it is to prove analogues of
Theorem~\ref{t:irrep} over fields of positive characteristic.
As too was already indicated (in~\S\ref{ss:iposchar}),  there is a
bonus to be had: our study enables a different approach to questions
about irreducibility of Specht modules.

As we observe in~\S\ref{ss:encode}, all
information about the action can conveniently be gathered together
into a matrix~$\bigmatrixl$ which breaks up nicely into blocks of the
same size (Proposition~\ref{p:bigmatrix}). The non-zero blocks all
lie along the diagonal and are all equal to a certain
matrix~$\matrixgl$ defined in~\S\ref{ss:cprod}. This matrix encodes
the multiplication table modulo lower cells of the~$C_w$ of
$\textup{\shape}~\lambda$. 
Theorem~\ref{t:formula} gives a formula for its determinant.

As remarked in~\S\ref{ss:ipf},  the argument in~\S\ref{ss:cprod} is
in effect deducing the
cellularity of the Kazhdan-Lusztig basis from results of~\cite{geck}
recalled in~\S\ref{s:setup};  and the argument in~\S\ref{ss:encode}
is in effect deducing Proposition~\ref{p:bigmatrix} as a consequence
of cellularity.
\mysubsection{On products of $C$-basis elements}\mylabel{ss:cprod}
Let $P_1$, \ldots, $P_m$ be the complete list of standard tableaux of shape~$\lambda$.    We claim:
\begin{multline}\label{e:cprod}
C(P_i,P_j)\cdot C(P_k,P_l)=
g_j^k C(P_i,P_l) \quad \bmod \\
\langle C_y\st 
\textup{\shape}(y)\triangleleft \lambda,\ 
y\klllneq (P_k,P_l),\   y\klrlneq (P_i,P_j)\rangle_A
\end{multline}  
the coefficient~$g_j^k$ being independent of $i$ and $l$.  

To prove the claim, consider the expression of the left hand side as
a linear combination of the $C$-basis elements.  For any~$C_y$ occurring
with non-zero coefficient, we have $y\klr (P_i,P_j)$ and $y\kll (P_k,P_l)$,
by the definition of the pre-orders (\S\ref{ss:klorders}).  
By~\S\ref{sss:kllr},
$\textup{\shape}(y)\trianglelefteq \lambda$;
and if $\textup{\shape}(y)\neq\lambda$,  then $y\klrlneq (P_i,P_j)$
and $y\klllneq (P_k,P_l)$.   If~$\textup{\shape}(y)=\lambda$,  then,
by~\S\ref{sss:cells}, $y\kllrsim(P_k,P_l)$;  by~\S\ref{sss:unrelate},
$y\kllsim(P_k,P_l)$; by~\S\ref{sss:cells}, the $Q$-symbol
of~$y$ is $P_l$;  and, analogously,  the $P$-symbol of~$y$ is~$P_i$.
That~$g_j^k$ doesn't depend upon~$i$ and~$l$ follows from the
description of the $\hecke$-isomorphisms between one sided cells
of the same \shape\ as recalled in~\S\ref{ss:cellmod},  and the claim
is proved.   We set 
\begin{equation}\label{e:matrixgl}
\matrixgl:=(g_j^k)_{1\leq j, k \leq m}
\end{equation}
\mysubsection{Relating the matrix $\matrixgl$ to the action on~$\rcmodl$}\mylabel{ss:encode}
%
Enumerate as $P_1$, \ldots, $P_m$ all the standard Young tableaux of
shape~$\lambda$.    Let us write $C(k,l)$ for the $C$-basis element
$C(P_k,P_l)$.  Consider the ordered basis 
$C(1,1)$, $C(1,2)$, \ldots, $C(1,m)$ of~$\rcmodl$.
Denote by~$\basise_i^j$ the element of~$\Endo{\rcmodl}$ that sends
$C(1,i)$ to $C(1,j)$ and kills the other basis elements.
Any element of~$\Endo{\rcmodl}$ can be written uniquely
as~$\sum \alpha_i^j\basise_i^j$.     Arrange the coefficients
as a row matrix like this:
\[
\left(
\begin{array}{cccc}
  \alpha_1^1 & \alpha_2^1 & \ldots & \alpha_m^1
\end{array} \st
\begin{array}{cccc}
  \alpha_1^2 & \alpha_2^2 & \ldots & \alpha_m^2
\end{array} \st
\begin{array}{ccc}
  \quad  & \ldots & \quad
\end{array} \st
\begin{array}{cccc}
  \alpha_1^m & \alpha_2^m & \ldots & \alpha_m^m
\end{array} 
\right)
\]

Now consider such row matrices for $\rhol(C(k,l))$.  Arrange
them one below the other,  the first row corresponding to the value
$(1,1)$ of $(k,l)$,  the second to $(2,1)$,  \ldots, the $m^{\textup{th}}$
row to $(m,1)$,  the ${(m+1)}^{\textup{th}}$ row to~$(1,2)$, \ldots,
and the last to $(m,m)$.   We thus get a matrix---denote it $\bigmatrixl$---of size
$\dl^2\times\dl^2$,  where $\dl:=\dim{\rcmodl}$.

Let us compute~$\bigmatrixl$ in the light of~(\ref{e:cprod}).
Setting $\alpha_i^j(k,l):=\alpha^j_i(\rhol(C(k,l)))$,  we have
(mind the abuse of notation:  this equation holds in~$\rcmodl$,
not in~$\hecke$):
\[ C(1,i)C(k,l)= \sum_j
\alpha_i^j(k,l)C(1,j).\]
Applying~(\ref{e:cprod}) to the left hand side and reading the result
as an equation in~$\rcmodl$,  we see that it equals $g_i^kC(1,l)$.   Thus
\[ \alpha_i^j(k,l)=\left\{
  \begin{array}{ll}
    g_i^k & \textup{if $j=l$}\\
    0 & \textup{otherwise}
  \end{array}
\right.\]
which means the following:
\begin{proposition}\mylabel{p:bigmatrix}
The matrix~$\bigmatrixl$ (defined earlier in this section)
is of block diagonal form,  with uniform
block size $\dl\times\dl$,  and each diagonal block equal to 
the matrix $\matrixgl=(g_i^k)$ of\/~\S\ref{ss:cprod}, where
the row index is~$k$ and the column index~$i$.
\end{proposition}

\mysubsection{A formula for the determinant of $\matrixgl$} %
\mylabel{ss:gdet}
\noindent 
Theorem~\ref{t:formula} below gives a formula for the determinant
of the matrix~$\matrixgl$ of~(\ref{e:cprod}).  
In order to state it, we need some notation.
Set
\begin{itemize}
\item
  $[\lambda]$:= the set of nodes in the Young diagram of shape $\lambda$;
\item
  $h_{ab}$:= hook length of the node $(a,b) \in [\lambda]$ (see~\S\ref{sss:dlambda}).
\item For a positive integer $m$,
\[\begin{split} \qintv{m}&:=v^{1-m} + v^{3-m}+ \cdots + v^{m-3} + v^{m-1}\\
\qintq{m}&:=1 + v^2 + v^4 + \cdots + v^{2(m-1)} \end{split}\]
\end{itemize}
Assuming $\lambda$ has $r$ rows, we can associate to~$\lambda$ a decreasing
sequence---called the {\em $\beta$-sequence\/}---of positive integers, 
the hook lengths
of the nodes in the first column of~$\lambda$.   The shape can
be recovered from the sequence, so the association gives a bijection
between shapes and decreasing sequences of positive integers.   
Given such a sequence $\beta_1>\ldots>\beta_r$,  write
$d(\beta_1,\ldots,\beta_r)$ for the number $d(\lambda)$  
of standard tableaux of shape~$\lambda$ (\S\ref{sss:dlambda}).  Extend
the definition of $d(\beta_1,\ldots,\beta_r)$ to an arbitrary sequence
of $\beta_1,\ldots,\beta_r$ of non-negative integers at most one of
which is zero as follows:
if the integers are not all distinct, then it is~$0$;    if the integers
are all distinct and positive,  then
it is $\sign(w)\,d(\beta_{w(1)},\ldots,\beta_{w(r)})$ where
$w$ is the permutation of the symmetric group~$\symm_r$
such that $\beta_{w(1)}>\ldots>\beta_{w(r)}$;    if the integers are distinct
and one of them---say~$\beta_k$---is zero,  then it is $d(\beta_1-1,\beta_2-1,
\ldots,\beta_{k-1}-1,\beta_{k+1}-1,\ldots, \beta_r-1)$,   which
is defined by induction on~$r$.
\begin{theorem}{\scshape (Hook Formula)}\mylabel{t:formula}  
For a partition~$\lambda$ of~$n$, 
\begin{equation}
  \label{eq:formula}
\det\matrixgl=  \epsilon_{w_{o, \lambda'}}^{d(\lambda)}
\prod\left(\frac{\qintv{h_{ac}}}{\qintv{h_{bc}}} \right)^{d(\beta_1,\ldots,
  \beta_a+h_{bc}, \ldots, \beta_b-h_{bc},\ldots, \beta_r)}
  \end{equation}
with notation as above,  
where
$\beta_1>\ldots>\beta_r$ is the $\beta$-sequence
of $\lambda$ and the product runs over 
$\{(a, b, c) \st (a,c), (b,c) \in [\lambda] \textrm{ and $a<b$}\}$.
\end{theorem}
\noindent
The proof of the theorem will be given in~\S\ref{s:gandgram} and
\S\ref{s:gdet}.   Some comments about it may be found at the
beginning of~\S\ref{s:gandgram}.

%% file: notation.tex
\mysection{Preliminaries about permutation and Specht modules}%
\mylabel{s:recall}\mylabel{s:notation}
\noindent
We recall some notation and results about permutation modules
and Specht modules needed in the sequel.
The isomorphism~$\theta$ recalled from~\cite{mp} in~\S\ref{ss:mpiso} plays a
fundamental role.

We keep the set up of~\S\ref{s:setup}.   We will also use freely the
notation and results recalled in~\S\ref{sss:invprelims}.   We will
be frequently referring to~\cite{kl,dj,mp}.  The reader
should be alert to the difference, pointed out in~\S\ref{sss:vandq}, 
between our notation and of these papers.

\mysubsection{Some notation}\mylabel{ss:notation}
Fix a partition~$\lambda\partition n$.   Let $\lambda'$, $\tupl$, $\tlowl$,
$\Wl$, $\wnotl$, $\dija$, $\xl$, 
and $\yl$ be as defined in~\S\ref{sss:invprelims}.
\begin{itemize}
\item
$\wl$ denotes the element of~$\dija$ that takes $\tupl$ to $\tlowl$.
By a {\em prefix\/} of~$\wl$ we mean an element of the 
form~$s_{i_1}\cdots s_{i_j}$ for some $j$, $1\leq j\leq k$,  where
$s_{i_1}\cdots s_{i_k}$ is some reduced expression of~$\wl$.
\item
Set $\zl:=v_{\wl} \xl\twl \ylp$.
\end{itemize}
It is elementary to see the following:
\begin{equation}
  \label{eq:wwl}
\textup{
$w_\lambda w$ belongs to $\dija$  and 
$\length(w_\lambda w)=\length(w_\lambda)+\length(w)$}\quad\quad
\textup{
for $w\in W_{\lambdap}$}.
\end{equation}
\mysubsection{Permutation modules~$\perml$ and Specht modules $\spechtl$}%
\mylabel{ss:permspecht}
Following~\cite{dj}---see \S3,~4 of that paper---we define
the {\em permutation module\/}~$\perml$ 
to be the right ideal~$\xl\hecke$,
the {\em Specht module\/}~$\spechtl$
to be the right ideal $\zl\hecke$.   They are $\hecke$-analogues
respectively of the $\symm_n$-representations on tabloids of
shape~$\lambda$ and its sub-representation the Specht module of
shape~$\lambda$ (defined respectively in~\S\ref{ss:tabloid} and
\S\ref{ss:specht}):   the corresponding $\symm_n$-modules
are recovered on setting~$v=1$ (see~\cite[Page~143]{mp}).
\mysubsubsection{Bases for~$\perml$ and~$\spechtl$}\mylabel{sss:mlslbase}
\begin{itemize}
\item \cite[Lemma~3.2~(1)]{dj} $\{\xl T_d\st d\in\dija\}$ is a basis for~$\perml$.
\item \cite[Theorem~5.6]{dj}
The elements $v_d\zl T_d$, $d$ a prefix of $\wlp$,  form an $A$-basis
for the Specht module~$\spechtl$ called the `standard basis'.
\end{itemize}
\mysubsubsection{The bilinear form $\langle\ ,\ \rangle$ 
on~$\perml$}\mylabel{sss:bilinear}
As in~\cite[page~34]{dj},  define a bilinear form 
$\langle\ \, ,\ \rangle$ on $\perml$ by setting
\(
\langle \xl T_d, \xl T_e\rangle
\)
equal to $1$ or~$0$ accordingly as elements $d$, $e$ of~$\dija$
are equal or not: as just recalled in~\S\ref{sss:mlslbase},
$\xl T_d$, $d\in\dija$,  form a basis for~$\perml$.
The form is evidently symmetric.  We have, by~\cite[Lemma~4.4]{dj}:
\begin{equation}\label{eq:bilinear}
\langle m_1h,m_2\rangle=\langle m_1,m_2h^*\rangle \quad
\textup{for
$m_1$, $m_2$ in~$\perml$ and $h$ in~$\hecke$}\end{equation} 
where $h\mapsto h^*$
is the $A$-anti-automorphism of~$\hecke$ given by $T_w\mapsto T_{w^{-1}}$.
\mysubsection{McDonough-Pallikaros isomorphism between right cell and Specht modules}\mylabel{ss:mpiso}
There is described in~\cite[Theorem~3.5]{mp}
a map from the right $\hecke$-ideal $\langle C_w\st w\klr \wnotl\rangle_A$
to $\permlp$:   it is denoted~$\theta$, defined by
$m\mapsto v_{\wlp}\xlp\twlp  m$,  and evidently a map 
of~$\hecke$-modules.%
\footnote{The factor~$v_{\wlp}$ in the definition of~$\theta$ appears
only in deference to~\cite{mp}.   If it
were omitted: the resulting~$\theta$ would still be a $\hecke$-module map;
the powers of~$v_{\wl}$ in Equations~(\ref{e:detgc})
and~(\ref{e:gandgram}) (and in the calculations leading up to them) would have to be omitted; and a suitable
power of~$v_{\wlp}$ would have to be added to Equation~(\ref{e:combine});
no other changes would have to be made.}
%
Using Proposition~\ref{p:pftinv2} and~(\ref{eq:xlyl}),
we can determine the image of~$\theta$:
$v_{\wlp}\xlp\twlp C_{\wnotl}\hecke=\xlp\twlp\yl\hecke=\spechtlp$.   
As to the kernel of~$\theta$,  it equals~$\langle C_w\st w\klrlt\wnotl\rangle_A$
as proved in~\cite[Theorem~3.5]{mp}.   Thus~$\theta$ gives an isomorphism
from the right cell module~$\rcmod(\wnotl)$ to $\spechtlp$.   Since~$\wnotl$ is
of shape~$\lambdap$ (Proposition~\ref{p:pftinv}~(2))
we conclude that $\rcmodl\cong\spechtl$.

%% file: gandgram.tex
\mysection{The first part of the proof of Theorem~\ref{t:formula}: relating $\det\matrixgl$ to the Gram determinant~$\gramdet$}%
\mylabel{s:gandgram}
\noindent
Towards the proof of Theorem~\ref{t:formula},
we relate,  using results from~\cite{dj,mp},  the determinant
of the matrix~$\matrixgl$ (defined in~(\ref{e:cprod})) to 
the {\em Gram determinant\/}~$\gramdet$,  namely, 
the determinant of the matrix of the restriction to the
Specht module~$\spechtl$ of the bilinear form~$\form$ 
on~$\perml$ (defined in~\S\ref{sss:bilinear}),  with respect
to the `standard basis' as in the second item in~\S\ref{sss:mlslbase}.   
The Gram determinant being well studied and results about it being
readily available in the literature,  we are thus lead 
to conclusions about~$\matrixgl$.
\ignore{ 
Towards establishing Theorem~\ref{t:formula}, we
relate the matrix $\matrixgl$ to the matrix of the ``Dipper-James
bilinear form'' on $\rcmodl$. By doing so, we are able to obtain a
relation between the determinant of $\matrixgl$ and the {\em Gram
  determinant\/}~$\gramdet$~--~the determinant of the matrix of the
restriction to the Specht module~$\spechtl$ of the bilinear
form~$\form$ on~$\perml$ (defined in~\S\ref{s:recall}), with respect
to the `standard basis' (ref.  Proposition\ref{p:recall}(12)). In
section \S\ref{ss:gtogram}, we give an explicit relation between the
two determinants: $\det\matrixgl$~and~$\gramdet$. A key step in obtaining this relation
comes from a result in \cite{mp}, using which we observe in
\S\ref{ss:tcbases} that the determinant of Dipper-James bilinear form on
$\rcmodl$ is the same with respect to its $C-$ and $T-$ bases.

We begin with some preliminary results and remarks. Here the reader
is reminded of the difference, as pointed out in~\S\ref{ss:hecke}, 
between our notation for the Hecke algebra and that of~\cite{kl,dj,mp}.
\input{tabprelim} 
}
\mysubsection{The Dipper-James bilinear form 
on~$\rcmodl$ computed in terms of its $C$-basis}%
\mylabel{ss:gramc}\noindent
Pulling back via the isomorphism~$\theta$ of~\S\ref{ss:mpiso} the
restriction to~$\spechtl$ of the bilinear form on~$\perml$ defined
in~\S\ref{s:notation},  we get a bilinear form on~$\rcmodl$ (which
we continue to denote by~$\form$).    Let us compute the matrix of this form
with respect to the basis $C(1,1)$, \ldots, $C(1,m)$,  where,
as in~\S\ref{ss:encode},  $P_1$, \ldots, $P_m$ is an enumeration of
all standard tableaux of shape~$\lambda$,  and $C(k,l)$ is short hand
notation for $C(P_k,P_l)$.
We further assume
that $P_1=t_\lambda$,  so that the right cell with $P$-symbol~$P_1$ is
the one containing $w_{0,\lambda'}$ (which under RSK corresponds to
the pair $(P_1,P_1)$---see Proposition~\ref{p:pftinv}~(2)).     The explanations for the steps in the
following calculation appear below:
\begin{multline*}
    \langle C(1,i), C(1,j)\rangle  = 
          \langle C(1,i)\theta, C(1,j)\theta\rangle
       = \langle \xl \twl \vwl C(1,i),\xl\twl \vwl C(1,j) \rangle\\
         = \vwl^2\langle \xl \twl, \xl\twl C(1,j) C(1,i)^* \rangle 
       = \vwl^2 \langle \xl \twl, \xl\twl C(1,j)C(i,1) \rangle \\
   = \vwl^2 \langle \xl \twl, \xl\twl g_j^iC(1,1) \rangle 
   = \vwl^2 g_j^i\, \langle \xl \twl, \xl\twl \epsilon_{\wnotlp}v_{\wnotlp}\ylp \rangle \\
   = \epsilon_{\wnotlp}v_{\wnotlp} \vwl^2 g_j^i\,\sum_{w\in W_{\lambdap}}\epsilon_w v_w^{-1}\langle \xl \twl, \xl \twl T_w \rangle 
  = \epsilon_{\wnotlp}v_{\wnotlp} \vwl^2 g_j^i
  \end{multline*}
The first equality follows from definition of the form on~$\rcmodl$; 
the second from the definition of~$\theta$;  the third from 
(\ref{eq:bilinear});
the fourth  
from (\ref{eq:c'star}).
For the fifth,
substitute for $C(1,j)C(i,1)$ using~(\ref{e:cprod}) and observe
that the `smaller terms' on the right hand side belong to the kernel
of~$\theta$ (\S\ref{ss:mpiso}).   The sixth follows by substituting for
$C(1,1)=C_{\wnotlp}$ from (\ref{eq:xlyl});
the seventh from the definition
of~$\ylp$;  and the final equality by combining the definition of the form
with (\ref{eq:wwl}) (observe that $\twl T_w=T_{\wl w}$ since $\length(\wl)
+\length(w)=\length(\wl w)$ and that $\wl w$ belongs to~$\dija$).

In particular,  the determinant of the matrix of the form~$\form$ on
$\rcmodl$ with respect to the basis $C(1,1)$, \ldots, $C(1,m)$ equals
\begin{equation}\label{e:detgc}
\epsilon_{\wnotlp}^{d(\lambda)}v_{\wnotlp}^{d(\lambda)} v_{\wl}^{2d(\lambda)}\det\matrixgl
\end{equation}
\mysubsection{The `$T$-basis' of $\rcmodlp$ and its relationship to the
$C$-basis}\mylabel{ss:tcbases}   Following~\cite[\S2]{mp},  we define the
`$T$-basis' of the right cell module~$\rcmodlp$ and show that it has
a uni-triangular relationship with the $C$-basis.     We do this by means
of the `$C$-basis' and `$T$-basis' of the right $\hecke$-module
$\cwnotl\hecke$,  which are defined respectively by:
\begin{itemize}
\item $C_w$, $w\klr\wnotl$ 
\item $\cwnotl T_d$, $d\in\dija$ 
\end{itemize}
That the `$C$-basis' is an $A$-basis follows from Proposition~\ref{p:pftinv2}.
That the `$T$-basis' is an $A$-basis is item~(\ref{i:tbasis}) in the
following:
\begin{proposition}\mylabel{p:srecall}
  \begin{enumerate}
\item\label{i:cwnotlty}
$\cwnotl T_y=\epsilon(y) v_y^{-1}\cwnotl$ for $y\in\Wl$.
\item\label{i:tbasis}
$\cwnotl T_d$, $d\in\dija$, form an $A$-basis for 
the right ideal $\cwnotl\hecke$.
\item \label{i:prefixdija}
$w\klrsim\wnotl$ if and only if $w=\wnotl d$ for a prefix $d$ of $\wl$.
In particular,  prefixes of~$\wl$ belong to $\dija$.
  \end{enumerate}
\end{proposition}
\begin{myproof}
(\ref{i:cwnotlty}) follows from \cite[Equation~(2.3.d)]{kl};
(\ref{i:tbasis}) from item~(\ref{i:cwnotlty}), Proposition~\ref{p:pftinv}~(\ref{i:dindija}),  
and~(\ref{e:ct})---see~\cite[Page~136]{mp};
(\ref{i:prefixdija}) from~\cite[Lemma~3.3~(iv)]{mp}.
\end{myproof}


The elements $w\klr\wnotl$ are precisely $\wnotl d$, $d\in\dija$
(Proposition~\ref{p:pftinv}~(\ref{i:dindija})).   Let $d_1$, \ldots, $d_M$
be the elements of~$\dija$ ordered so that $i\leq j$ 
if $d_i\leq d_j$ in the Bruhat order.  By~\cite[Proposition~2.13]{mp}
and its proof,   the two bases above are related by a uni-triangular
matrix with respect to an ordering as above (keeping in mind that
$T_y$ in the notation of~\cite{mp} equals $v_yT_y$ in ours):
\[
\left(
  \begin{array}{c}
    \cwnotl T_{d_1}\\
\vdots\\
\cwnotl T_{d_M}
  \end{array}
\right)
           =  \left(  
             \begin{array}{ccc}
               1 &  & 0\\
                 & \ddots & \\
               \star & & 1 
             \end{array}
\right)
\left(
  \begin{array}{c}
    C_{\wnotl d_1}\\
\vdots\\
C_{\wnotl d_M}
  \end{array}
\right)
\]

Let us now read this equation in the quotient~$\rcmodlp$ of~$C_{\wnotl}\hecke$.
Let $d_{i_1}$, \ldots, $d_{i_m}$ with $1\leq i_1<\ldots<i_m\leq M$
be such that they are all the prefixes of~$\wl$---see 
Proposition~\ref{p:srecall}~(\ref{i:prefixdija})---so that 
$\wnotl d_{i_1}$, \ldots, $\wnotl d_{i_m}$ are all the
elements right equivalent to~$\wnotl$.
Writing $e_1$, \ldots, $e_m$
in place of $d_{i_1}$, \ldots, $d_{i_m}$,    and noting that
$C_{\wnotl d_j}$ vanishes in~$\rcmodlp$ unless $\wnotl d_j\klrsim\wnotl$, we have:
\[
\left(
  \begin{array}{c}
    \cwnotl T_{e_1}\\ 
\vdots\\
\cwnotl T_{e_m}
  \end{array}
\right)
           =  \left(  
             \begin{array}{ccc}
               1 &  & 0\\
                 & \ddots & \\
               \star & & 1 
             \end{array}
\right)
\left(
  \begin{array}{c}
    C_{\wnotl e_1}\\ 
\vdots\\
C_{\wnotl e_m}
  \end{array}
\right)
\]

We conclude that $\cwnotl T_{e_1}$, \ldots, $\cwnotl T_{e_m}$ form an
$A$-basis for $\rcmodlp$.  It is called the {\em $T$-basis\/} and is
in uni-triangular relationship with the $C$-basis $C_{\wnotl e_1}$,
\ldots, $C_{\wnotl e_m}$.   In particular, the determinants of the
matrices of the form~$\form$ on~$\rcmodlp$ (defined in~\S\ref{ss:gramc}) 
with respect to the $T$- and $C$-bases are the same.
\mysubsection{$\det\matrixgl$ and the Gram
determinant~$\gramdet$}\mylabel{ss:gtogram}
Continuing towards our goal of relating $\det\matrixgl$ to 
the Gram determinant~$\gramdet$,   let us compute the image under the map~$\theta$
of the $T$-basis elements of~$\rcmodl$.     Given a prefix~$e$ of
$\wlp$,   we have 
\[
\begin{array}{rcll}
C_{\wnotlp}T_e\theta & =  & v_{\wl}\xl \twl C_{\wnotlp} T_e & 
\textrm{(by the definition of~$\theta$ in~\S\ref{ss:mpiso})}\\
& = &
\epsilon_{\wnotlp} v_{\wnotlp} \left(v_{\wl} \xl \twl \ylp\right) T_e 
& \textrm{(by (\ref{eq:xlyl}))}\\
& =  & \epsilon_{\wnotlp} v_{\wnotlp} \zl T_e 
& \textrm{(by the definition of $\zl$ in \S\ref{ss:notation})}\\
& =  & \epsilon_{\wnotlp} v_{\wnotlp} v_e^{-1}(v_e \zl T_e)\\
\end{array} 
\]
Noting that $v_e \zl T_e$ is a standard basis element of~$\spechtl$
(see \S\ref{sss:mlslbase}), we conclude that
the determinant of the matrix of the bilinear form~$\form$ on~$\rcmodl$
with respect to the $T$-basis equals 
$v_{\wnotlp}^{2d(\lambda)} 
\left(\prod_{e} v_e^{-1}\right)^2 \gramdet$.
Combined with~(\ref{e:detgc}) and the conclusion
of~\S\ref{ss:tcbases},   this gives
\begin{equation}\label{e:gandgram}
\det\matrixgl=
(\epsilon_{\wnotlp}
v_{\wnotlp}  v_{\wl}^{-2})^{d(\lambda)}
(\prod_{e} v_e)^{-2} \gramdet
\end{equation}
where the product is taken over all prefixes~$e$ of~$\wlp$.

%% file: gdet.tex
\mysection{Conclusion of the proof of Theorem~\ref{t:formula}} %
\mylabel{s:gdet}
\noindent    In this section,  we complete the proof of 
Theorem~\ref{t:formula} by using results from 
\cite{djblocks,jmathas}.   Both sides
of Equation~(\ref{eq:formula}) are elements of~$A$.
To prove they are equal,  we may pass
to the quotient field $K:=\rational(v)$ of~$A$.
We do this tacitly in the sequel.    Observe that $\hecke_K$ is semisimple
(see~\ref{sss:hspecial}).
\mysubsection{Second half of the proof of Theorem~\ref{t:formula}}%
\mylabel{s:2ndhalf}
We set things up to be able to use a
formula from~\cite{djblocks} for the Gram  determinant $\gramdet$.
Let $S_1$, \ldots,
  $S_m$ be an enumeration of all the standard tableaux of
  shape~$\lambda$.  For $i$, $u$ such that $1\leq i\leq m$, $1\leq
  u\leq n$, let $S_i^u$ denote the standard tableau obtained
  from~$S_i$ by deleting all nodes with entries exceeding~$u$;
  set $\gamma_{ui}:=\prod_{j=1}^{a-1}
  {\qintq{h_{jb}}}/{\qintq{h_{jb}-1}}$ where $(a,b)$ is the position
  of the node in~$S_{i}^u$ containing~$u$, $h_{jb}$ is the hook
  length in~$S_i^u$ of the node in
  position~$(j,b)$, and 
  $\qintq{s}:=1+v^2+v^4+\cdots+v^{2(s-1)}$ for a positive integer~$s$.
  By~\cite[Theorem~4.11]{djblocks}, the Gram determinant~$\gramdet$ is
  given by
\begin{equation}
  \label{e:djblocks}
  \gramdet = v^{2r} \prod_{i=1}^m \prod_{u=1}^n \gamma_{ui}
\quad\quad\quad\textrm{for some integer $r$}
\end{equation}

We now apply the equation in~\cite[Corollary~2.30, page~251]{jmathas}.
Computing~$\Delta_\mu(\lambda')$ (in the notation of~\cite{jmathas})
with~$\mu=1^n$, and re-indexing the product in the right side of that
equation over nodes of~$\lambda$ rather than of~$\lambda'$, we get
\begin{equation}\label{e:jmathas}
  \prod_{i=1}^m \prod_{u=1}^n \gamma_{ui} =
\prod\left(\frac{\qintq{h_{ac}}}{\qintq{h_{bc}}} \right)^{d(\beta_1,\ldots,
  \beta_a+h_{bc}, \ldots, \beta_b-h_{bc},\ldots, \beta_r)}
\end{equation}
where 
the product on the right hand side 
runs over triples $(a,b,c)$ as in the statement of the theorem.

Combining Equations~(\ref{e:gandgram}), (\ref{e:djblocks}), (\ref{e:jmathas}),
and (\ref{e:r}),   we get 
\begin{equation}
  \label{e:combine}
\epsilon_{\wnotlp}^{d(\lambda)} v_{\wnotlp}^{d(\lambda)}\det\matrixgl
=
\prod\left(\frac{\qintq{h_{ac}}}{\qintq{h_{bc}}} \right)^{d(\beta_1,\ldots,
  \beta_a+h_{bc}, \ldots, \beta_b-h_{bc},\ldots, \beta_r)}
\end{equation}
The left hand side is an element of~$A$.  
As to the right hand side,
it is regular with value~$1$ at $v=0$,
since the same is true for $\qintq{s}$ for every positive integer~$s$.
Thus both sides of the equation 
belong to~$1+v\mathbb{Z}[v]$ and \[\det\matrixgl=
\epsilon_{\wnotlp}^{d(\lambda)} v_{\wnotlp}^{-d(\lambda)}+\textrm{higher degree terms}.\]
The `bar-invariance' of the $C$-basis elements (\S\ref{ss:hecke})
means that:
\[\textrm{$\bar{g_j^k}=g_j^k$ for $g_j^k$ as in~(\ref{e:cprod}) and
so also $\bar{\det\matrixgl}=\det\matrixgl$.}\]
Thus $\det\matrixgl$ has the form: 
\[ \epsilon_{\wnotlp}^{d(\lambda)} v_{\wnotlp}^{-d(\lambda)}+\cdots+
\epsilon_{\wnotlp}^{d(\lambda)} v_{\wnotlp}^{d(\lambda)}\] 
the terms represented by $\cdots$ being of $v$-degree strictly between
$-d(\lambda)\length(\wnotlp)$ and $d(\lambda)\length(\wnotlp)$.  Equating the
$v$-degrees on both sides of~(\ref{e:combine}) gives
\begin{equation*}
  d(\lambda)\length(\wnotlp)=
\sum
d(\beta_1,\ldots,
  \beta_a+h_{bc}, \ldots, \beta_b-h_{bc},\ldots, \beta_r)
\left(h_{ac}-h_{bc}\right).
\end{equation*}
Using this and substituting~$v^{h_{ac}}\qintv{h_{ac}}$, $v^{h_{bc}}\qintv{h_{bc}}$, 
respectively for $\qintq{h_{ac}}$, $\qintq{h_{bc}}$ into (\ref{e:combine}),
we arrive at the theorem.
  \begin{flushright}
$\Box    $
  \end{flushright}
\begin{lemma}
  \mylabel{l:r}
The integer~$r$ in the exponent of~$v$ in Equation~(\ref{e:djblocks})
is given by
\begin{equation}
  \label{e:r}
  r = d(\lambda)\left(\length(\wl)-\length(\wnotlp)\right)+
\sum_{i=1}^m \length(d_i)
\end{equation}
where $d_1$, \ldots, $d_m$ are all the prefixes of~$\wlp$.
\end{lemma}
\begin{myproof}
  We essentially work through the proof
  of~\cite[Theorem~4.11]{djblocks} to calculate the exponent of $v$
  appearing in \ref{e:djblocks}.  

  Let $d_1$, \ldots, $d_m$ 
  be ordered so that $i<j$ if $\length(d_i)<\length(d_j)$.  Let $e_i$
  and $f_i$, $1\leq i\leq m$, be bases of~$\spechtl$ as
  in~\cite{djblocks}.  The $e_i:=v_{d_i}\zl T_{d_i}$ are just the
  standard basis (see~\S\ref{sss:mlslbase}).  The $f_i$ are an
  orthogonal basis in uni-triangular relationship with
  the~$e_i$~\cite[Theorem~4.7]{djblocks}.  Thus
  $\gramdet=\prod_{i=1}^{m} \langle f_i, f_i\rangle$.

  Let the enumeration $S_1$, \ldots, $S_m$ of standard tableaux of
  shape~$\lambda$ be such that $S_i=t_\lambda d_i$.  Then,
  from~\cite[Lemma~4.10]{djblocks}, $\langle f_i, f_i\rangle
  =v^{2r_i}\prod_{u=1}^{n}\gamma_{ui}$.  We claim:
  \begin{enumerate}
  \item $r_1=\length(\wl)-\length(\wnotlp)$
  \item for $i>1$, $r_i=r_j+1$ where $j<i$ such that $e_i=v e_j
    T_{(k-1,k)}$.
  \end{enumerate}
  The lemma being clear given the claim, it remains only to prove the
  claim.

  Item~(2) of the claim follows from the following two observations
  made in the course of the proof of~\cite[Lemma~4.10]{djblocks}:
  $\langle f_i, f_i\rangle= c_j \langle f_j, f_j \rangle$, and
  \[
  \gamma_{ui}=\left\{
    \begin{array}{cl}
      \gamma_{uj} & \textrm{if $u\neq k-1, k$}\\
      \gamma_{kj} & \textrm{if $u=k-1$}\\
      v^{-2}c_j\gamma_{k-1,j} & \textrm{if $u=k$}\\
    \end{array}
  \right.
  \]
  (The definition of~$c_j$ is irrelevant for our purposes.)

  To prove item~(1) of the claim, we compute $\langle f_1,
  f_1\rangle$.
  We have $f_1=e_1=\zl$.  Substituting for $\zl$ and in turn for $\yl$
  from their definitions in~\S\ref{ss:notation} and \S\ref{sss:invprelims}, 
we get
  \[
  \begin{array}{rcl}
    \langle f_1, f_1 \rangle 
    & = & \langle \zl, \zl \rangle\\
    & = & \langle v_{\wl}\xl T_{\wl}\ylp, v_{\wl}\xl T_{\wl}\ylp\rangle\\
    & = & v_{\wl}^2\sum_{u,u'\in W_{\lambdap}}\epsilon_u\epsilon_{u'}v_u^{-1} v_{u'}^{-1} \langle \xl T_{\wl}T_u, \xl T_{\wl}T_{u'}\rangle\\
  \end{array} \]
  Using in order~(\ref{eq:bilinear}),~(\ref{eq:tutu'}),~(\ref{eq:wwl}) 
  and the definition of~$\form$,  we get:
  \[
  \begin{array}{rcl}
    \langle \xl T_{\wl}T_u, \xl T_{\wl}T_{u'}\rangle & = &
    \langle \xl T_{\wl}, \xl T_{\wl}T_{u'}T_{u^{-1}}\rangle\\
    & = & \langle \xl T_{\wl}, \xl T_{\wl}(T_{u'u^{-1}}+\sum_{w\in W_{\lambdap}; w>u'u^{-1}}c_w T_w)\rangle\\
    & = & \langle \xl T_{\wl}, \xl T_{\wl u'u^{-1}}+\sum_{w\in W_{\lambdap}; w>u'u^{-1}}c_w T_{\wl w}\rangle\\
    & = & \left\{ 
      \begin{array}{cl}
        1 & \textrm{if $u=u'$}\\
        0 & \textrm{otherwise}
      \end{array}
    \right.
  \end{array}
  \]
  so that \[ \langle f_1, f_1 \rangle = \vwl^2 \sum_{u\in
    W_{\lambdap}}v_u^{-2} =\vwl^2 \vwnotlp^{-2} \sum_{u\in
    W_{\lambdap}}\vwnotlp^{2}v_u^{-2} =\vwl^2 \vwnotlp^{-2} \sum_{u\in
    W_{\lambdap}}v_u^{2}\] Routine calculations show:
\[\sum_{u\in \symmn}v_u^{2}\ \ =\ \
\qintqfac{n}\quad\quad\quad\textrm{ and }\quad\quad\quad \sum_{u\in
  W_{\lambdap}}v_u^{2}\ \ =\ \
\qintqfac{\lambdap_1}\cdots\qintqfac{\lambdap_r}
\]
where $\qintqfac{n}:= \qintq{n}\qintq{n-1}\cdots\qintq{1}$ and
$\lambda'=(\lambdap_1,\ldots,\lambdap_r)$.  Finally, a pleasant
verification, given the fact that $S_1=t_\lambda$, shows:
\[
\qintqfac{\lambdap_1}\cdots\qintqfac{\lambdap_r} = \prod_{u=1}^n\gamma_{u1}
\]
The proof of the claim (and so also of the lemma) is complete.
\end{myproof}

%% file: poschar.tex
\mysection{On the irreducibility of Specht modules}%
\mylabel{s:poschar}
\noindent
Let $\ring$ be a field and $a$ a non-zero element of~$\ring$. Fix
notation as in~\S\ref{s:setup}, and consider~$\ring$ as an $A$-module via
the map $A\to\ring$ defined by $v\mapsto a$.  

The purpose of this
section is to revisit the question of when the Specht module
$\spechtlk:=\spechtl\tensor_A\ring$ (which, by~\S\ref{ss:mpiso}, is
isomorphic to the right cell module $\rcmodl\tensor_A\ring)$ is
irreducible (as a module over
$\hecke_\ring:=\hecke\tensor_Ak$).  As we will see, 
the results of~\S\ref{s:matrixgl}--\ref{s:gdet} 
afford us a fresh approach to this question.
In addition, they allow us to prove
generalizations of Theorems~\ref{p:klbasis},~\ref{p:p2kl} to situations
when $\heckek$ is not necessarily semisimple.

The form $\form$ defined on~$\perml$ by Dipper-James has been
recalled in~\S\ref{sss:bilinear}.   Let $\form_\ring$ denote the
form on~$\perml\tensor_A\ring$ obtained by extension of scalars and also
its restriction to $\spechtlk$.
Note that $\form_\ring$ is symmetric (since $\form$ on $\perml$ is).
\mysubsection{Preliminaries}
\mylabel{ss:prelim}
We first establish some notation and make a few observations.
Thus equipped, we recall some results from~\cite{murphy95} 
in a form that is convenient or us.

For an automorphism~$\ddag$ of~$\hecke$ and 
a (right) $\hecke$-module~$M$, we denote by $M^\ddag$ the (right) 
$\hecke$-module whose underlying $A$-module is the same as that 
of~$M$---it is convenient
to write $m^\ddag$ for an element~$m$ of~$M$ thought of 
in~$M^\ddag$---and with the action
of~$\hecke$ being given by $m^\ddag h^\ddag:=(mh)^\ddag$.   
\begin{itemize}
\item 
For a (right) $\hecke$-ideal~$\ideali$, the symbol~$\ideali^\ddag$
can be interpreted without conflict as either the image of $\ideali$
under~$\ddag$ or the module~$M^\ddag$ defined as above taking $M$ to
be~$\ideali$.
\item 
If~$\ddag$ is an involution,  then $M\simeq (M^\ddag)^\ddag$ naturally.
\item 
For an anti-automorphism~$\ddag$ of~$\hecke$,  $M^\ddag$ defined similarly would
naturally be a left~$\hecke$-module:   $h^\ddag m^\ddag:=(mh)^\ddag$.
\end{itemize}

Now consider the involution~$\dag$ on~$\hecke$ defined in \S\ref{sss:2inv1aa}
and the permutation module~$\perml$ defined in~\S\ref{ss:permspecht}.
From equations~(\ref{eq:cc'char}),~(\ref{eq:cjc'}), and~(\ref{eq:xlyl}),
it follows that  $(C_{\wnotl}\hecke)^\hash=\perml$;
from Proposition~\ref{p:pftinv2} that
\begin{equation}\label{eq:permlbase}
\textup{$C_w'$, $w\klr\wnotl$, form an $A$-basis for~$\perml$.}
\end{equation}
Set~$\nlambda:=\langle C_w'\st \textup{\shape}(w)\domleq\lambdap\rangle_A$
and~$\nhatl:=\langle C_w'\st \textup{\shape}(w)\doml\lambdap\rangle_A$.
From~\S\ref{sss:kllr} and Proposition~\ref{p:pftinv}~(\ref{i:wnotl})
it follows that $\perml\subseteq\nlambda$.
Set~$\sbl:=\perml/\perml\cap\nhatl$.   
From~\S\S\ref{sss:cells}--\ref{sss:unrelate} it follows that
$\{w\st w\klr\wnotl\textup{ and }\textup{\shape}(w)\doml\lambdap\}
=\{w\st w\klrlneq\wnotl\}$, so that $\perml\cap\nhatl
=\langle C_w'\st w\klrlneq\wnotl\rangle_A$. Thus the images
in $\sbl$ of $C_w'$, $w\klrsim\wnotl$, form a basis for~$\sbl$.
And we get
\begin{equation}
  \label{eq:sandsb}
\rcmodlp^\hash\simeq\sbl \quad\quad\quad\quad\quad \rcmodl \simeq (\sblp)^\hash
\end{equation}

For a (right) $\hecke$-module~$M$,  the dual $M^\dual:=\Hom_A(M,A)$ is
naturally a left $\hecke$-module:  $(m)(h\phi):=(mh)\phi$, for 
$\phi\in M^\dual$, $m\in M$, and $h\in\hecke$.   We use the
anti-automorphism~$*$ defined in~\S\ref{sss:2inv1aa} to switch between
right and left module structures: $(M^\dual)^*$ becomes a right $\hecke$-module.
\begin{itemize}
\item 
The process $M\mapsto M^\dual$ commutes with that of~$M\mapsto M^\ddag$ defined
earlier in this section:  in particular, $(M^\dual)^*\simeq(M^*)^\dual$ naturally.
\item
If~$M$ is free as an $A$-module,  then $(M^\dual)^\dual\simeq M$ naturally.
\end{itemize}
\begin{proposition}\textup{(\cite[Theorem~5.2]{murphy95})}
  \mylabel{p:dual}
We have an isomorphism $((\spechtl)^\dual)^*\simeq(\spechtlp)^\hash$.  
In particular $\spechtlk$
is irreducible if and only if~$\spechtlpk$ is so.
\end{proposition}
\begin{myproof}
It is proved in~\cite[Theorem~5.2]{murphy95} that
$(\sblp)^\hash\simeq((\sbl)^\dual)^*$:  a perfect pairing
$(\ ,\,\,): \sbl\times(\sblp)^\hash\to A$ with the property that
$(m,nh)=(mh^*,n)$ is given.  
Combining this statement with
the isomorphisms~(\ref{eq:sandsb}) and
the isomorphism $\rcmodl\simeq\spechtl$
of~\S\ref{ss:mpiso}, the proposition follows.
\end{myproof}

A shape~$\lambda$ is called {\em $e$-regular\/} if the number of rows in it
of any given length is less than~$e$.
Let now $e$ be the smallest positive 
integer such that $1 + a^2 + \cdots + a^{2(e-1)} =0$;  
if there is no such integer,  then $e=\infty$.  
\begin{proposition}\textup{(\cite[Theorem 6.9]{murphy95})}
  \mylabel{p:nonzero}
If~$\lambda$ is $e$-regular,  the bilinear form~$\form_k$
on~$\spechtlk$ is non-zero.
\end{proposition}
\begin{myproof}
  Consider the form~$\forml$ on $\perml$ defined in~\cite[page~114]{murphy95}.
Using (\ref{eq:c'star}),  this can be expressed
in our notation as follows in terms of the basis
(\ref{eq:permlbase}) of $\perml$:   for~$w$, $x$ such that
$w\klr\wnotl$, $x\klr\wnotl$, $\langle C_w',C_x'\rangle_\lambda$ is the
coefficient of~$C_{\wnotl}'$ in $C_w'C'_{x^{-1}}$. 
It follows readily from the definition
in~\S\ref{ss:klorders} of the relation~$\klr$ that
if either $w\klrlneq\wnotl$ or $x\klrlneq\wnotl$ (which is equivalent
to $x^{-1}\klllneq\wnotl^{-1}=\wnotl$),  then $\langle C'_w,C_x'\rangle_\lambda=0$.
Thus $\forml$ descends to~$\sbl$.

From~\cite[Theorem 6.9---see also its proof]{murphy95} it follows that
$\forml$ does not vanish on $\sbl$ if $\lambda'$ is $e$-regular.
This means that there exist $w\klrsim\wnotl$, 
$x\klrsim\wnotl$ such that
the coefficient of $C_{\wnotl}'$ in $C_w'C_{x^{-1}}'$ is non-zero.
Applying the involution~$\hash$,  we conclude that the coefficient of
$C_{\wnotl}$ in $C_wC_{x^{-1}}$ is non-zero (see~(\ref{eq:c'dag})).   

Consider the ordered pairs of standard tableaux associated to 
$\wnotl$, $w$, and $x$ by the RSK-correspondence:  $\wnotl\leftrightarrow
(\tlowlp,\tlowlp)$ by Proposition~\ref{p:pftinv}~(\ref{i:wnotl});
let~$Q_w$ and $Q_x$ be the standard tableaux of shape~$\lambdap$ such
that $w\leftrightarrow(\tlowlp,Q_w)$ and $x\leftrightarrow(\tlowlp,Q_x)$
(see~\S\ref{ss:cellsandkrs}).   Then $x^{-1}\leftrightarrow(Q_x,\tlowlp)$
(see~\S\ref{ss:rsk}) and, by (\ref{e:cprod}), 
\[C_wC_{x^{-1}}=C(\tlowlp,Q_w)C(Q_x,\tlowlp)\equiv g_w^x C(\tlowlp,\tlowlp)
\bmod\textup{``lower terms''}\]
The conclusion of the last paragraph translated to this notation says
that the coefficient $g_w^x$ is non-zero.

Consider the pull-back to $\rcmodlp_\field$ via the isomorphism~$\theta$
of~\S\ref{ss:mpiso} of the form $\langle\ ,\ \rangle_\field$ on~$\spechtlpk$.
Denoting it too by $\langle\ ,\ \rangle_\field$,  
the big display in~\S\ref{ss:gramc} says that 
$\langle C_w,C_x\rangle_\ring$ equals the coefficient~$g_w^x$ up to
sign and a power of~$v$.   Thus $\langle C_w,C_x\rangle_\field\neq0$,
which means that the form $\langle\ ,\ \rangle_\field$ on~$\spechtlpk$
is non-zero.\end{myproof}
\mysubsection{Analogues of
Theorems~\ref{p:klbasis},~\ref{p:p2kl} for not necessarily 
semi-simple~$\heckek$}

\begin{theorem}\mylabel{t:cpklbasis} 
For an $e$-regular shape~$\lambda$ such that $\spechtlk$ is
irreducible, 
the Kazhdan-Lusztig basis elements $C_w$, $w$ of \shape~$\lambda$, 
thought of as operators on~$\spechtlk$ form a basis for $\Endo{\spechtlk}$.
\end{theorem}
\begin{myproof}
  By (\ref{eq:bilinear}), the radical of the form $\form_\ring$ on
$\spechtlk$ is a
  $\heckek$-submodule.  Since $\spechtlk$ is assumed irreducible, the form
is  either identically zero or non-degenerate.  But, as shown in
  Proposition~\ref{p:nonzero} above, it is non-zero under the
  assumption of $e$-regularity of~$\lambda$.  Thus its matrix with
  respect to any basis of~$\spechtlk$ has non-zero determinant.  By
  (\ref{e:detgc}), $\det\matrixgl\vert_{v=a}$ is such a determinant
  (up to a sign and power of~$a$), so it is non-zero.  It now follows
  from Proposition~\ref{p:bigmatrix} that the operators~$C_w$, $w$ of
  \shape~$\lambda$, form a basis for~$\Endo{\spechtlk}$.\end{myproof}
\begin{corollary}
  \mylabel{c:cpklbasistwo} Suppose that $\lambdap$ is $e$-regular and
  that~$\spechtlk$ is irreducible.  Then the elements $C_w'$,
  $\textup{\shape}(w)=\lambdap$, as operators on~$\spechtlk$ form a
  basis for $\Endo{\spechtlk}$.
\end{corollary}
\begin{myproof}
  By Theorem~\ref{t:cpklbasis},  the $C_w$, $\textup{\shape}(w)=\lambdap$,
as operators on~$\spechtlpk$ form a basis for $\Endo{\spechtlpk}$ ($\spechtlpk$
is irreducible by Proposition~\ref{p:dual}).  
Since $\spechtlpk\simeq ((({\spechtlk})^\dual)^*)^\hash$ 
(Proposition~\ref{p:dual} again),  and ${C_w}^\dag=\epsilon_w C_w'$ by
(\ref{eq:c'dag}),
it follows that the $C_w'$, $\textup{\shape}(w)=\lambdap$, as operators
on~$((\spechtlk)^\dual)^*$ form a basis for $\Endo{((\spechtlk)^\dual)^*}$.
Since $(C_w')^*=C'_{w^{-1}}$ by (\ref{eq:c'star}) and the \shape{s} of
$w$ and $w^{-1}$ are the same,  the result follows.
\end{myproof}
\begin{theorem}\mylabel{t:cpp2kl}
  Let~$\subsets$ be the set of $e$-regular shapes~$\lambda$ such that
  the Specht module~$\spechtlk$ is irreducible.  Let $U$ be a finite
  dimensional semisimple $\hecke_\ring$-module, every irreducible
  component of which is of the form~$\spechtlk$, $\lambda\in\subsets$.
  Let $\subsett$ be the subset of~$\subsets$ consisting of those
  shapes~$\lambda$ such that~$\spechtlk$ appears as a component
  of~$U$.  Then the images in $\Endo{U}$ of\/ $C_x$, $x\in\symm_n$
  such that\/ $\textup{\shape}(x)$ belongs to $\subsett$, form a basis
  for the image of~$\hecke_\ring$ in~$\Endo{U}$ (under the map
  $\hecke_\ring\to\Endo{U}$ defining~$U$).
\end{theorem}
\begin{myproof}
  The proof is similar to that of Theorem~\ref{p:p2kl}.
\end{myproof}

\mysubsection{A criterion for irreducibility of
  $\spechtlk$} \mylabel{ss:spechtirred}
We first observe 
that Proposition~\ref{p:bigmatrix} gives us a criterion
for irreducibility of Specht modules (Theorem~\ref{t:cpirred1}).
We then deduce from this criterion a conjecture of Carter
\cite[Conjecture~1.2]{james:carter} about irreducibility of Specht modules
(Corollary~\ref{t:cpirred}).
Of course the conjecture has long been proved~\cite{james:carter,jmathas},
but our approach is new.
\begin{theorem}\mylabel{t:cpirred1}  If $\det\matrixgl|_{v=a}$ does
not vanish in~$k$,  then $\spechtlk$ is irreducible.  
\end{theorem}
\begin{myproof}
  Suppose that $\det\matrixgl|_{v=a}$ does not vanish in~$\ring$.
Then, by Proposition~\ref{p:bigmatrix},  the matrix~$\bigmatrix$
is invertible (in~$\ring$, after specializing to~$v=a$).
Thus the elements $C_w$, $w$ of \shape~$\lambda$, are linearly
independent (and so form a basis) as operators on~$\spechtlk$.
In particular,  $\spechtlk$ is irreducible,  and the assertion is proved.
\end{myproof}
\mysubsubsection{A new proof of Carter's conjecture}\mylabel{sss:carter}
Let $p$ denote the smallest positive integer such that $p=0$ in~$\ring$;
if no such integer exists, then $p=\infty$.
For an integer~$h$,  define $\nu_p(h)$ as the 
largest power of~$p$ (possibly~$0$) that divides~$h$ in case~$p$ is
positive,  and as~$0$ otherwise.     
Recall that $e$ denotes the smallest positive 
integer such that $1 + a^2 + \cdots + a^{2(e-1)} =0$;  
if there is no such integer,  then $e=\infty$.  For an integer~$h$, define
\[
\nu_{e,p}(h):=\left\{
  \begin{array}{ll}
    0 & \textrm{if $e=\infty$ or $e\nmid h$}\\
    1+\nu_p(h/e)  & \textrm{otherwise}\\
  \end{array}
\right.
\]
The {\em $(e,p)$-power diagram\/} of shape~$\lambda$ is the 
filling up of the nodes of the shape~$\lambda$
by the $\nu_{e,p}$'s of the respective hook lengths. 

Observe that $e=p$ if $a=1$.   
\begin{corollary}\mylabel{t:cpirred} 
\textup{\cite{james:carter,jmathas}}
If 
the $(e,p)$-power
diagram of~$\lambda$ has either no column or no row
containing different numbers,  
then $\spechtlk$ is irreducible.   
\end{corollary}
\begin{myproof}
It is enough to do the case when 
no column of the $(e,p)$-power diagram has different numbers:  if the condition
is met on rows and not on columns,  we can pass to~$\lambda'$ and use
the observation (\cite[Corollary~3.3]{fl} or Proposition~\ref{p:dual} above) 
that $\spechtl_k$ is irreducible if and only if $\spechtlp_k$ is. 

So assume that in every column of the $(e,p)$-power diagram 
the numbers are all the same.   
We claim that each of the factors $\qintv{h_{ac}}/\qintv{h_{bc}}$ 
on the right hand side 
of~(\ref{eq:formula}) makes sense as an element of~$k$ and 
is non-zero.   Combining the claim with 
Theorems~\ref{t:formula} and~\ref{t:cpirred1}
yields the assertion.

To prove the claim,  we need 
the following elementary observations,  where $h$ denotes a positive
integer:
\begin{itemize}
\item $\qintv{h}$ vanishes in~$\ring$ if and only if $e$ is finite and 
divides $h$.
\item
if $e$ is finite and divides~$h$,  then $\qintv{h}=(\qintv{h/e})\mid_{v=v^e}\qintv{e}$.   
\item 
$a^{2e}=1$ if $e$ is finite.
\end{itemize}
If either $e=\infty$ or $e$ does not divide any of the hook lengths
in shape~$\lambda$,  then the claim follows from the
the first of the above observations.
So now suppose that~$e$ is finite and 
divides either $h_{ac}$ or $h_{bc}$.    By our hypothesis, $e$ then divides
both $h_{ac}$ and $h_{bc}$;    moreover both $h_{ac}/e$
and $h_{bc}/e$ are divisible by~$p$ to the same extent.   Using the second
and third observations above,  we conclude that the image in~$\ring$ of
$\qintv{h_{ac}}/\qintv{h_{bc}}$ is the same as that of the 
rational number $h_{ac}/h_{bc}$ (written in reduced form), and so is non-zero.   
\end{myproof}

%% file: bibliography.tex

	\newcommand\citenumfont[1]{\textbf{#1}}

\bibliographystyle{bibsty-final-no-issn-isbn}
\addcontentsline{toc}{section}{References}
\ifthenelse{\equal{\finalized}{no}}{
\bibliography{abbrev,references}
}{\input{bbl-file.tex}}